\documentclass[preprint,3p,12pt]{elsarticle}
\biboptions{numbers,sort&compress}
\usepackage{bm}
\usepackage{amssymb}
\usepackage{amsthm,amsmath}
\usepackage{hyperref}
\usepackage{booktabs}
\usepackage[hyphenbreaks]{breakurl}
\usepackage{color}
\usepackage{graphicx}
\usepackage{subfigure}
\usepackage{rotating}

\numberwithin{equation}{section}

\newtheorem{theorem}{Theorem}[section]
\newtheorem{proposition}{Proposition}[section]
\newtheorem{lemma}{Lemma}[section]
\newtheorem{remark}{Remark}[section]
\newtheorem{corollary}{Corollary}[section]

\journal{Mathematical Methods in Applied Sciences}
\begin{document}
\begin{frontmatter}



\title{Fast implicit difference schemes for time-space fractional diffusion
equations with the integral fractional Laplacian}

\author[a,b,c]{Xian-Ming Gu}
\ead{guxianming@live.cn, guxm@swufe.edu.cn}
\author[b]{Hai-Wei Sun}
\cortext[cor1]{Corresponding author}
\ead{hsun@umac.mo}
\author[d]{Yanzhi Zhang}
\ead{zhangyanz@mst.edu}
\author[e]{Yong-Liang Zhao\corref{cor1}}
\ead{ylzhaofde@sina.com}
\address[a]{School of Economic Mathematics/Institute of Mathematics,\\
Southwestern University of Finance and Economics, Chengdu, Sichuan 611130, China}
\address[b]{Department of Mathematics, University of Macau, Avenida da Universidade,
Taipa, Macao, China}
\address[c]{Bernoulli Institute for Mathematics, Computer Science and Artificial Intelligence,\\
 University of Groningen, Nijenborgh 9, P.O. Box 407, 9700 AK Groningen, The Netherlands}
\address[d]{Department of Mathematics and Statistics,\\
 Missouri University of Science and Technology, Rolla, MO 65409-0020, United States}
\address[e]{School of Mathematical Sciences,\\ University of Electronic Science and
Technology of China, Chengdu, Sichuan 611731, China}
\begin{abstract}
In this paper, we develop two fast implicit difference schemes for solving a class of variable-coefficient
time-space fractional diffusion equations with integral fractional Laplacian (IFL). The proposed schemes
utilize the graded $L1$ formula for the Caputo fractional derivative and a special finite difference
discretization for IFL, where the graded mesh can capture the model problem with a weak singularity
at initial time. The stability and convergence are rigorously proved via the $M$-matrix
analysis, which is from the spatial discretized matrix of IFL. Moreover, the proposed schemes use the
fast sum-of-exponential approximation and Toeplitz matrix algorithms to reduce the computational
cost for the nonlocal property of time and space fractional derivatives, respectively. The fast schemes greatly reduce the computational work of solving
the discretized linear systems from $\mathcal{O}(MN^3 + M^2N)$ by a direct solver to $\mathcal{O}(MN(\log N + N_{exp}))$ per
preconditioned Krylov subspace iteration and a memory requirement from $O(MN^2)$ to $O(NN_{exp})$, where $N$ and $(N_{exp} \ll)~M$
are the number of spatial and temporal grid nodes. The spectrum of preconditioned matrix is also given for ensuring
the acceleration benefit of circulant preconditioners. Finally, numerical results are presented to show the utility of the proposed
methods.
%
%
\end{abstract}
\end{frontmatter}


\section{Introduction}
\label{sec1}
In recent decades, fractional partial differential equations (FPDEs) have attracted growing attention in
modeling phenomena with long-term memory and spatial heterogeneity arising in engineering,
physics, chemistry and other applied sciences \cite{Podlub,Sun18}. In physics, fractional derivatives are used
to model anomalous diffusion. Anomalous diffusion is the theory of diffusing particles in environments that
are not locally homogeneous \cite{zquez,Li2019x,Chen2010,MFSBJ}. A physical-mathematical model to anomalous
diffusion may be based on FPDEs containing derivatives of fractional order in both space and time, where the
sub-diffusion appears in time and the super-diffusion occurs in space simultaneously \cite{Rodriguez,Chi2017}. On the other
hand, although most of time-space fractional diffusion models are initially defined with the spatially integral
fractional Laplacian (IFL) \cite{Du2012,Tian16,Duo2019,Teso18}, many previous studies (cf. e.g., \cite{Podlu09,Chen2010,Yang10x,Yang11,Li2019x}) always
substitute the space Riesz fractional derivative \cite{Podlub} for the IFL. In fact, such
two kinds of definitions are
not equivalent in high-dimensional cases \cite{Yang11,Garbac,Duo2019}. It means that the `direct' study of time-space fractional
diffusion models with the IFL should be worthily considered.
%

In this paper, we study an alternative time-space fractional diffusion equation (TSFDE) with variable
coefficients in one space dimension
\begin{equation}
\begin{cases}
{}^{C}_0D^{\gamma}_t u(x,t) = -\kappa(x,t)(-\Delta)^{\alpha/2}u(x,t) + f(x,t),&x\in\Omega,~~t\in(0,T],\\
u(x,t) = 0,&x\in\Omega^c,~~t\in[0,T],\\
u(x,0) = \phi(x),& x\in\Omega
\end{cases}
\label{eq1.1}
\end{equation}
where $\kappa(x,t) > 0$ denotes the diffusivity coefficients,
$\Omega = (-l,l)$, $\Omega^{c} = \mathbb{R}\setminus\Omega$, and the initial condition $\phi(x)$ and the source term $f(x,t)$ are known functions. Meanwhile,
${}^{C}_0D^{\gamma}_t$ is the Caputo derivative \cite{Podlub} of order $\gamma$ - i.e.
\begin{equation}
{}^{C}_0D^{\gamma}_tu(x,t) =
\begin{cases} \frac{1}{\Gamma(\gamma-1)}\int^{t}_0\frac{1}{(t - s)^{\gamma}}\frac{\partial
u(x,s)}{\partial s}ds,& 0 < \gamma < 1,\\
\frac{\partial u(x,t)}{\partial t}, &\gamma = 1,
\label{eq1.2xx}
\end{cases}
\end{equation}
and throughout the paper we always assume that $0 < \gamma < 1$. Here the fractional Laplacian $(-\Delta)^{\alpha/2}$ is defined by \cite{Guan05,Du2012,Tian16,Duo2019}:
\begin{equation}
(-\Delta)^{\alpha/2}u(x) = c_{1,\alpha}~{\rm P.V.}\int_{\mathbb{R}}\frac{u(x) - u(x')}{|x - x'|^{1 + \alpha}}
dx',\quad \alpha\in(0,2),
\label{eq1.4}
\end{equation}
where P.V. stands for the Cauchy principal value,  and $|x - x'|$ denotes  the  Euclidean  distance between
points $x$ and $x'$.  The normalization constant $c_{1,\alpha}$ is defined as
\begin{equation}
c_{1,\alpha} =
\frac{2^{\alpha - 1}\alpha\Gamma(\frac{\alpha + 1}{2})}{\sqrt \pi\Gamma(1 - \alpha/2)}
\end{equation}
with $\Gamma(\cdot)$ denoting the Gamma function. From a probabilistic point of view, the IFL represents the infinitesimal generator of a symmetric $\alpha$-stable L\'{e}vy process
\cite{Guan06,Duo2019,Garbac}. Mathematically, the well-posedness/regularity of the Cauchy problem or uniqueness
of the solutions of the TSFDE (\ref{eq1.1}) has been studied in \cite{Koloko,zquez,Hanyga,Li2018x,Padgett,Jia2016}.

Due to the nonlocality, the analytical (or closed-form) solutions of TSFDEs (\ref{eq1.1}) on a finite domain
are rarely available. Therefore, we have to rely on numerical treatments that produce approximations to the
desired solutions; refer e.g., to \cite{Podlub,Bonito17x,Teso18,Bonito18,Acosta} and references therein for a description of such approaches. In
fact, utilizing the suitable temporal discretization, most of the early established numerical methods including
the finite difference (FD) method \cite{Arshad17,Arshad,Chi2017,Gu2017}, finite element (FE) method \cite{Yue2018,Feng16}, and
matrix (named it as all-at-once) method \cite{Podlu09,Biala18} for the TSFDE (\ref{eq1.1}) were developed
via the fact that the IFL is equivalent to the Riesz fractional derivative in one space dimension \cite{Yang10x}. However,
such a numerical framework cannot be directly extended to solve the two- and three-dimensional TSFDEs due
to the IFL $(-\Delta)^{\alpha/2}u(x,y)\neq -\frac{\partial^{\alpha}u(x,y)}{\partial |x|^{\alpha}} - \frac{\partial^{
\alpha}u(x,y)}{\partial |y|^{\alpha}}$ \cite{Yang11,Garbac,Duo2019}. Therefore, it will hinder the development of numerical solutions
for TSFDEs from the stated objective.

In order to remedy the above drawback, Duo, Ju and Zhang \cite{Duo2018b} replace the IFL in TSFDE (\ref{eq1.1})
by the spectral fractional Laplacian \cite{Yang11,Duo2019} and present a fast numerical approach which combines
the matrix transfer method \cite{Yang11} with inverse Laplace transform for solving the one- and multi-dimensional TSFDEs
(\ref{eq1.1}) with \textit{constant coefficients}. Although the numerical results show that their proposed
method converges with the second-order accuracy in both time and space variables, the spectral fractional
Laplacian on a bounded domain is also not equivalent to the IFL at all \cite{Duo2019}. On the other hand,
Nochetto, Ot\'{a}rola and Salgado \cite{Nochetto} use the Caffarelli-Silvestre extension to rewrite the
TSFDE (\ref{eq1.1}) with $\kappa(x,t)\equiv \kappa$ as a two-dimensional quasi-stationary elliptic problem
with dynamic boundary condition. Then, they establish a FE scheme for solving the converted elliptic problem
and show that the numerical scheme cannot reach the error estimates of order $O(\tau^{2-\gamma})$ claimed in
the literature. Later, Hu, Li and Li \cite{Hu2017,Hu2018} successively exploit the similar strategy with FD approximation
for the converted elliptic problem of one- and multi-dimensional TSFDEs (\ref{eq1.1}) with $\kappa(x,t)
\equiv \kappa$. Nevertheless, the numerical results show that such FD schemes often converge with the
less than first- and second-order accuracy in time and space, respectively, even for TSFDEs with sufficient
smooth solutions.

In fact, it is important to set up numerical schemes which utilize the `direct' discrezations of IFL for
solving the TSFDEs (\ref{eq1.1}). Moreover, the discretizations of (multi-dimensional) become a recently hot
topic, with the main numerical challenge stemming from the approximation the hypersingular integral, see e.g.
\cite{Huang14,Tian16,Duo2018a,Duo2019z,Hao2019,Bonito19,Acosta3,Acosta4,Minden,Ainsworth}. Indeed, there are some numerical
schemes utilized the temporal $L1$ formula \cite{Lin07} (or numerical Laplace inversion \cite{Bonito17x}) and
spatial FE discretization \cite{Tian16,Bonito18,Bonito19,Acosta3,Acosta4} for solving the (multi-dimensional)
constant-coefficient TSFDEs (\ref{eq1.1}) \cite{Bonito18,Rodriguez,Acosta}. Both the theoretical and numerical
results are reported to show that such numerical schemes are efficient to solve the (multi-dimensional) TSFDEs
(\ref{eq1.1}) with $\kappa(x,t)\equiv \kappa$. In addition, there are some other kinds of time-space fractional
diffusion models but related to TSFDEs (\ref{eq1.1}), where the spatial (or temporal) nonlocal operator is a
replacement for the IFL (or the Caputo fractional derivative). This is mainly because the nonlocal operators
with suitable kernels can exactly embrace the IFL and the Caputo fractional derivative, respectively \cite{Du2012,Guan15,Tian16,Liu2017,Chen17}. For such novel
model problems, Guan and Gunzburger \cite{Guan15} establish a class of numerical methods unitized the $\theta$ schemes
and piecewise-linear FE discretization. Their fully discrete scheme is analyzed for all to determine conditional
and unconditional stability regimes for the scheme and also to obtain error estimates for the approximate solution.
Later, Liu, et al. \cite{Liu2017} improve the idea of Guan and Gunzburger by giving the proof of convergence behavior with
$\mathcal{O}(\tau^{2 - \gamma} + h^2)$. Meanwhile, Liu, et al. consider the piecewise-quadratic FE discretization
to improve the spatial convergence rate. The efficient implementation based on fast Toeplitz-matrix multiplications
\cite{Chan96,Ng2004,Lei2013} of their proposed scheme is also reported. For space-time nonlocal diffusion equations,
Chen, et al. \cite{Chen17} propose a numerical scheme, by exploiting the quadrature-based FD method in time and the Fourier
spectral method in space, and show its stability. Moreover, it is shown that the convergence is uniform at a rate
of $\mathcal{O}(\delta + \sigma^2)$ (where $\delta$ and $\sigma$ are the time and space horizon parameters) under
certain regularity assumptions on initial and source data. Even these are several methods with linear solvers of
quasilinear complexity, the implementation of the above methods is still complicated, especially the computation
of entries of stiffness matrix in FE discretization or finding the modes in terms of expansion basis in spectral
method, cf. \cite{Liu2017,Chen17,Bonito18,Minden}. In particular, it is pointed out that ``More than 95\% of CPU time is used to assembly routine" for their
FE methods \cite{Acosta4}.

On the other hand, most of the above mentioned methods overlook that the presence of the kernel $(t - s)^{-\gamma}
$ results in a weak initial singularity in the solution of Eq. (\ref{eq1.1}), so that approximation methods (e.g.,
$L1$ formula) on \textit{uniform meshes} have a poor convergent rate and high computational cost \cite{Sakamoto,Nochetto,Chen17,Shen18}. In this work, we devote us to developing fast implicit difference schemes (IDSs) for solving the TSFDE (\ref{eq1.1}), the direct scheme
utilizes the simple FD discretization \cite{Duo2018a} and the graded $L1$ formula \cite{Shen18} for approximating
the IFL and the Caputo fractional derivative, where the non-uniform temporal discretization can overcome the initial
singularity. Due to the repeated summation of numerical solutions in the previous steps, the direct scheme always needs
much CPU time and memory cost, especially for the larger number of time steps \cite{Fu2017a,Fu2017b,Gu2020}. In order to alleviate the computational
cost, the sum-of-exponential (SOE) approximation \cite{Jiang17} of the kernel $(t - s)^{-\gamma}$ in the graded $L1$ formula for Caputo
fractional derivative can be efficiently evaluated via the recurrence method. Thus, we can derive the fast implicit
difference scheme. In particular, we revisit the matrix properties of the discretized IFL and prove the discretized
matrix is a strictly diagonally dominant and symmetric $M$-matrix with positive diagonal elements (i.e., the symmetric
positive definite matrix), which is not studied in \cite{Duo2018a}. Based on such matrix properties, we strictly prove
that the fast numerical schemes for the TSFDE (\ref{eq1.1}) are unconditionally stable and present the corresponding
error estimates of $\mathcal{O}(M^{-\min\{r\gamma,2-\gamma\}} + h^2)$ ($h$ is the spatial grid size) under certain regularity assumptions on the
smooth solutions. To our best knowledge, there are few successful attempts to derive the efficient IDSs for solving the TSFDE (\ref{eq1.1})
with rigorous theoretical analyses. This is one of main attractive advantages of our proposed methods compared to the above mentioned methods.

In addition, the nonlocality of IFL results in dense discretized linear systems, which is the leading time-consuming
part in practical implementations \cite{Bonito18,Acosta4,Minden}. Fortunately, the coefficient matrix of discretized linear
systems enjoys the Toeplitz-like structure \cite{Duo2018a,Duo2019z,Hao2019}, it means that we solve the sequence of discretized
linear systems in a matrix-free pattern \cite{Pang12,Lei2013,Gu2017}, because the Toeplitz
matrix-vector products can be computed via fast Fourier transforms (FFTs) in $\mathcal{O}(N\log N)$ operations. More precisely, we
will adapt the circulant preconditioners \cite{Chan96,Ng2004,Lei2013} for accelerating the Krylov subspace solvers \cite{Ng2004} for the sequence of discretized
linear systems. Moreover, the benefit of circulant preconditioners will be verified via both theoretical and numerical
results. It notes that fast schemes greatly reduce the computational work of solving the discretized linear systems from $\mathcal{
O}(MN^3 + M^2N)$ by a direct solver to $\mathcal{O}(MN(\log N + N_{exp}))$ per preconditioned Krylov subspace iteration and a memory
requirement from $O(MN^2)$ to $O(NN_{exp})$ (see Section \ref{sec3.1} for defining $N_{exp}$).

The contributions of the current work can be summarized as follows.
\begin{itemize}
  \item We present two IDSs for solving the TSFDE (\ref{eq1.1}) with non-smooth initial data and such numerical
schemes can be easily extended to solve the multi-dimensional cases.
  \item Both the stability and convergence of these IDSs are rigorously proved via the discrtized
matrix properties.
  \item We provide the efficient implementation of fast IDSs with theoretical guarantee for reducing
  the computation and memory cost deeply.
\end{itemize}

The rest of this paper is organized as follows. In Section \ref{sec3}, both direct and fast IDSs are
derived for the TSFDE (\ref{eq1.1}) in details, and their stability and convergence are proved by revisiting the properties
of spatial discretized matrix. In Section \ref{sec4}, the efficient implementation based on fast preconditioned
Krylov subspace solvers of the proposed IDSs are given and the accelerating benefit of circulant preconditioners
is theoretically guaranteed via clustering eigenvalues around 1. Section \ref{sec5} presents numerical results
to support our theoretical findings and show the effectiveness of the proposed IDSs. Finally, the paper closes with
conclusions in Section \ref{sec6}.
%
%
%
%
\section{Direct and fast implicit difference schemes}
\label{sec3}
In this section, we will establish two implicit difference schemes for solving the problem (\ref{eq1.1}).
Meanwhile, the stability, convergence and error analysis of such difference schemes are investigated and
proved in details.
\subsection{Two implicit difference schemes}
\label{sec3.1}
As mentioned above, we assume that the problem (\ref{eq1.1}) has a solution $u(x,t)$ such that
\begin{equation}
\left|\frac{\partial^k u(x,t)}{\partial t^k}\right|\leq c_0t^{\gamma - k},\quad 0\leq
k\leq 3,
\label{eq2.1}
\end{equation}
here and in what follows, $\hat{C}$ and $c_j~(j=0,1,2)$ are positive constants,
which depend on the problem but not on the mesh parameters \cite{Shen18}. Let $M,
N,r\in\mathbb{N}^{+}$ (positive integers), $h = 2l/N, x_i = -l + ih$, $t_m = (m/M)^rT$ and $\tau_m
 = t_m - t_{m - 1}$. We also consider the sets
\begin{equation*}
\Omega_h = \{x_i|0\leq i\leq N\},~~ \Omega_{\tau} = \{t_m|0\leq m\leq M\},~~
\Omega_{h,\tau} = \{(x_i,t_m)|0\leq i\leq N,~0\leq m\leq M\}
\end{equation*}
and let ${\bm v} = \{v^{m}_i| 0 \leq i \leq N,~0 \leq m \leq M\}$ be a grid function on $\Omega_{
h,\tau}$. We define the set $\mathcal{V}_h = \{{\bm v}|{\bm v} = (v_0, v_1,\cdots, v_{
N-1}, v_N)\}$ of grid functions on $\Omega_h$ and provide it with the norm
\begin{equation*}
\|{\bm v}\|_{\infty} = \max_{0\leq i\leq N}|v_i|.
\end{equation*}

For $m\geq 1$, we approximate the Caputo fractional derivative (\ref{eq1.2xx}) by the $L1$
formula on the graded mesh, which can capture the weak initial singularity of (\ref{eq1.1}):
\begin{equation*}
{}^{C}_0D^{\gamma}_tu(t_m) \approx D^{\gamma}_tu(t_m) \triangleq \frac{1}{\Gamma(1 - \gamma)}\Big[a^{(m,\gamma)
}_mu(t_m)- \sum^{m-1}_{k = 1}(a^{(m,\gamma)}_{k+1} - a^{(m,\gamma)}_k)u(t_k) - a^{(m,\gamma)}_1u(t_0)
\Big],
\end{equation*}
where
\begin{equation}
a^{(m,\gamma)}_k = \frac{1}{\tau_k}\int^{t_k}_{t_{k-1}}\frac{ds}{(t_m - s)^{\gamma}},~~k\geq 1.
\label{eq2.8}
\end{equation}
The truncation error $\psi^{m}$ can be defined by $\psi^{m}\triangleq {}^{C}_0D^{\gamma}_tu(t_m) - D^{\gamma}_tu(t_m)$.
From \cite[Lemma 2.1]{Shen18}, we obtain the boundness of the truncation error
\begin{equation}
|\psi^{m}|\leq \hat{C}m^{-\min\{r(1 + \gamma),2 - \gamma\}},
\quad m = 1,2,\cdots,M,
\label{lem2.1}
\end{equation}
if $|u''(t)|\leq c_0t^{\gamma-2},~0\leq t\leq T$ and $\hat{C} > 0$.

On the other hand, it notes that the above graded $L1$ scheme always
needs much computational cost in practical applications due to the repeatedly weighted sum of the solutions of
previous time steps. To reduce the cost, here it is useful to develop the fast approximation
of Caputo fractional derivative on a non-uniform temporal mesh.
\begin{lemma}{\rm(\cite{Jiang17})}
Let $\epsilon, \delta$ and $T$ denote the tolerance error, cut-off time restriction and final
time, respectively. Then there exist $N_{exp}\in\mathbb{N}^{+}$ and $s_j, w_j > 0,~j
= 1, 2,\cdots,N_{exp}$ such that
\begin{equation*}
\left|t^{-\gamma} - \sum^{N_{exp}}_{j=1}w_je^{-s_jt}\right| \leq \epsilon,\quad{\rm for~
any}~ t\in[\delta,T],
\end{equation*}
where $N_{exp} = \mathcal{O}((\log\epsilon^{-1})(\log\log\epsilon^{-1} + \log(T\delta^{-1
})) + (\log\delta^{-1})(\log\log\epsilon^{-1} + \log\delta^{-1}))$.
\label{lem2.2}
\end{lemma}

Based on Lemma \ref{lem2.2}, we set $\delta = (1/M)^r T$, and then the fast approximation of
Caputo fractional derivative on a graded temporal grid can be drawn as follows ($m\geq 1$).
\begin{equation}
\begin{split}
{}^{C}_0D^{\gamma}_tu(t_m) & = {}^{FC}D^{\gamma}_tu(t_m) + \mathcal{O}(m^{-\min\{r(1 + \gamma),2-\gamma\}} + \epsilon)\\
& \triangleq \frac{1}{\Gamma(1 - \gamma)}\Big[b^{(m,\gamma)}_mu(t_m) -\sum^{m-1}_{k=1}(b^{(m,\gamma)}_{k+1}
- b^{(m,\gamma)}_k)u(t_k) - b^{(m,\gamma)}_1u(t_0)\Big]\\
&~~~+~\mathcal{O}(m^{-\min\{r(1 + \gamma),2-\gamma\}} + \epsilon),\quad m = 1,2,\cdots,M,\\
\end{split}
\label{lem2.3}
\end{equation}
where the estimate of truncation error holds when $|u'(t)|\leq c_0t^{\gamma-1}$ and $|u''(t)|\leq c_0t^{\gamma-2}$ and
\begin{equation}
b^{(m,\gamma)}_k =
\begin{cases}
\sum\limits^{N_{exp}}_{j = 1}w_j\frac{1}{\tau_k}\int^{t_k}_{t_{k-1}}e^{-s_j(t_m - s)}ds,&
k = 1,2,\cdots,m-1,\\
a^{(m,\gamma)}_m,& k = m.
\end{cases}
\end{equation}

Moreover, we provide the information for approximating the IFL with extended
Dirichlet boundary conditions in (\ref{eq1.1}). According to the idea in \cite{Duo2018a},
the approximation is given by
\vspace{-1mm}
\begin{equation}
\begin{split}
(-\Delta)^{\alpha/2}_{h,\mu}u(x_i) & = C^{h}_{\alpha,\mu}\Bigg[\Big(\sum\limits^{N-1}_{\ell
= 2}\frac{(\ell + 1)^{\nu} - (\ell - 1)^{\nu}}{\ell^{\mu}} + \frac{N^{\nu} - (N-1)^{\nu}}{
N^{\mu}} + (2^{\nu} + \kappa_{\mu} - 1) \\
&~~~+ \frac{2\nu}{\alpha N^{\alpha}}\Big)u(x_i) - \frac{2^{\nu} + \kappa_{\mu} - 1}{2}\Big(
u(x_{i-1}) + u(x_{i+1})\Big)\\
&~~~- \frac{1}{2}\sum^{N-1}_{j=1,j\neq i,i\pm 1}\frac{(|j-i| + 1)^{\nu} - (|j-
i|-1)^{\nu}}{|j - i|^{\mu}}u(x_j)\Bigg],\\
\end{split}
\label{eq2.11}
\end{equation}
where $i = 1,2,\cdots,N-1$, $C^{h}_{\alpha,\mu} = c_{1,\alpha}/(\nu h^{\alpha}) > 0$ and the constant
$\kappa_{\mu} = 1$ for $\mu\in(\alpha,2)$, while $\kappa_{\mu} = 2$ if $\mu = 2$. Meanwhile, we denote
$\nu = \mu - \alpha$ for notational simplicity.
\begin{lemma}{\rm(\cite{Duo2018a})}
Suppose that $u(x)\in\mathcal{C}^{s,\frac{\alpha}{2}}(\mathbb{R})$ has the finite support on an open set
$\Omega\subset\mathbb{R}$, and Eq. (\ref{eq2.11}) is a finite difference approximation of the
fractional Laplacian $(-\Delta)^{\alpha/2}$. Then, for any $\mu\in(\alpha,2)$, there is
\begin{equation}
\|(-\Delta)^{\alpha/2}u(x) - (-\Delta)^{\alpha/2}_{h,\mu}u(x)\|_{\infty,\Omega}\leq \tilde{C}h^p
\end{equation}
with $\tilde{C} > 0$ depending on $\alpha$ and $\mu$.  Here $p\in(0,2]$ would be
determined via the regularity (i.e., the index $s\in\mathbb{N}$) of $u(x)$.
\label{lem2.4}
\end{lemma}
\vspace{-2mm}
Lemma \ref{lem2.4} provides a direct discretization for the IFL appeared in the TSFDE (\ref{eq1.1}).
At present, the spatial and temporal discretizations are ready for developing the numerical methods.
Evaluating Eq. (\ref{eq1.1}) at the points $(x_i, t_m)$, we have
\begin{equation}
{}^{C}_0D^{\gamma}_tu(x_i,t_m) = -\kappa(x_i,t_m)(-\Delta)^{\alpha/2}u(x_i,t_m) + f(x_i,t_m),
\end{equation}
where $1\leq i\leq N-1, 1 \leq m\leq M$. Let ${\bm u}=\{u^{m}_i\}_{ 0 \leq i \leq N,0\leq m\leq M}$ be a grid function defined by
\begin{equation*}
U^{m}_i:= u(x_i,t_m),\quad f^{i}_m = f(x_i,t_m),\quad \kappa^{m}_i = \kappa(x_i,t_m),\quad 0\leq i\leq N, 0 \leq m\leq M.
\end{equation*}
Using these notations and recalling Eq. (\ref{lem2.3}) along with Lemma \ref{lem2.4}, we can
approximate Eq. (\ref{eq1.1}) at grid point $(x_i, t_m)$ as follows
\begin{equation}
\begin{cases}
{}^{FC}D^{\gamma}_t U^{m}_i = -\kappa^{m}_i(-\Delta)^{\alpha/2}_{h,\mu}U^{m}_i +  f^{m}_i +
R^{m}_i,&1\leq i\leq N-1,~1\leq m\leq M,\\
U^{m}_0 = U^{m}_{N} \equiv 0,& 0\leq m \leq M,\\
U^{0}_i = \phi(x_i),& 1\leq i\leq N - 1,
\end{cases}
\label{eq3.14}
\end{equation}
where the terms $\{R^{m}_i\}$ are small and satisfy the inequality
\begin{equation*}
|R^{m}_i|\leq c_2(m^{-\min\{r(1 + \gamma),2-\gamma\}} + h^p + \epsilon),
\quad 1 \leq i\leq N-1,~1\leq m\leq M.
\end{equation*}
We omit the above small terms and arrive at the following implicit difference scheme
\begin{equation}
\begin{cases}
\frac{1}{\Gamma(1 - \gamma)}\left[b^{(m,\gamma)}_mu^{m}_i -
\sum\limits^{m-1}_{k=1}(b^{(m,\gamma)}_{k+1} - b^{(m,\gamma)}_k)u^{k}_i - b^{(m,\gamma)}_1
u^{0}_i\right] = -\kappa^{m}_i(-\Delta)^{\alpha/2}_{h,\mu}u^{m}_i \\
\qquad\qquad\qquad\qquad\qquad\qquad\qquad\qquad\quad~+~f^{m}_i,\quad 1\leq i\leq N-1,~1\leq m\leq M,\\
u^{m}_0 = u^{m}_{N} = 0,\qquad\qquad\qquad\qquad\qquad\qquad\qquad\quad~~0\leq m \leq M,\\
u^{0}_i = \phi(x_i),\qquad\qquad\qquad\qquad\qquad\qquad\qquad\qquad~~1\leq i\leq N - 1,
\end{cases}
\label{eq3.15}
\end{equation}
which is named as the fast implicit difference scheme (FIDS). Similarly, we combine the graded $L1$ formula with Lemma \ref{lem2.4} for deriving the following difference scheme
\begin{equation}
\begin{cases}
\frac{1}{\Gamma(1 - \gamma)}\left[a^{(m,\gamma)}_mu^{m}_i -
\sum\limits^{m-1}_{k=1}(a^{(m,\gamma)}_{k+1} - a^{(m,\gamma)}_k)u^{k}_i - a^{(m,\gamma)}_1
u^{0}_i\right] = -\kappa^{m}_i(-\Delta)^{\alpha/2}_{h,\mu}u^{m}_i \\
\qquad\qquad\qquad\qquad\qquad\qquad\qquad\qquad\quad~+~f^{m}_i,\quad 1\leq i\leq N-1,~1\leq m\leq M,\\
u^{m}_0 = u^{m}_{N} = 0,\qquad\qquad\qquad\qquad\qquad\qquad\qquad\quad~~0\leq m \leq M,\\
u^{0}_i = \phi(x_i),\qquad\qquad\qquad\qquad\qquad\qquad\qquad\qquad~~1\leq i\leq N - 1,
\end{cases}
\label{eq3.16}
\end{equation}
which is labelled as direct implicit difference scheme (DIDS).

At each time level, both FIDS (\ref{eq3.15}) and DIDS (\ref{eq3.16}) are the resultant linear
systems, which can be solved by the direct method (e.g., Gauss elimination) with total computational
cost $\mathcal{O}(N^3M + NMN_{exp})$ for FIDS and $\mathcal{O}(N^3M + NM^2)$ for DIDS. Note that,
generally, $N_{exp} < 100$ \cite{Jiang17,Shen18} and (if) $M$ is very large, so that FIDS requires smaller
computational cost than DIDS. Moreover, FIDS only requires $O(N^2 + NN_{exp})$ memory units rather than $\mathcal{
O}(N^2 + NM)$ for DIDS. In Section \ref{sec4}, we will further reduce the computational
cost of both FIDS and DIDS by means of matrix-free preconditioned iterative solvers.

\subsection{The stability and convergence}
\label{sec3.2}
In this subsection, we discuss the stability and convergence of the difference scheme for
the problem (\ref{eq1.1}). In order to analyze the stability and convergence, we rewrite
the FIDS (\ref{eq3.15}) into the matrix form
\begin{equation}
\mathcal{M}^{(m)}{\bm u}^{m} = \frac{1}{\Gamma(1 - \gamma)} \left[\sum\limits^{m-1}_{k=1}
(b^{(m,\gamma)}_{k+1} - b^{(m,\gamma)}_k){\bm u}^{k} + b^{(m,\gamma)}_1{\bm u}^{0}\right]
+ {\bm f}^{m},
\label{eq3.16gu}
\end{equation}
where $\mathcal{M}^{(m)} = \frac{1}{\Gamma(1 - \gamma)}b^{(m,\gamma)}_mI + K^{(m)}A$ and refer to
Eq. (\ref{eq3.16x}) for the definition of $A$, $I$ is the
identity matrix of order $N - 1$, $K^{(m)} = {\rm diag}(\kappa^{m}_1,\kappa^{m}_2,\cdots,\kappa^{
m}_{N-1})$, ${\bm u}^m = [u^{m}_1,u^{m}_2,\cdots,u^{m}_{N-1}]^\top$, ${\bm f}^m = [f^{m}_1,f^{m}_2,
\cdots,f^{m}_{N-1}]^\top$, and $b^{(m,\gamma)}_m > 0$ \cite[Lemma 2.4]{Shen18}. First of all, we revisit the properties of spatial discretization,
which is not deeply studied in the original paper \cite{Duo2018a,Duo2019z}. In fact, the spatial discretization
(\ref{eq2.11}) of $(-\Delta)^{\alpha/2}u(x,t)$ can be expressed in the matrix-vector product form $(-
\Delta)^{\alpha/2}_{h,\mu}{\bm u}^m = A{\bm u}^m$, where $A = [a_{ij}]_{i,j=1,\cdots,N-1}$ is the matrix representation of the (discretized) fractional
Laplacian, defined as
\begin{equation}
a_{ij} = C^{h}_{\alpha,\mu}\begin{cases}
\sum\limits^{N-1}_{\ell = 2}\frac{(\ell + 1)^{\nu} - (\ell - 1)^{\nu}}{\ell^{\mu}} + \frac{N^{\nu}
- (N-1)^{\nu}}{N^{\mu}} + (2^{\nu} + \kappa_{\mu} - 1) + \frac{2\nu}{\alpha N^{\alpha}},& j = i,\\
- \frac{2^{\nu} + \kappa_{\mu} - 1}{2}, & j = i\pm 1,\\
-\frac{(|j - i| + 1)^{\nu} - (|j - i| - 1)^{\nu}}{2|j - i|^{\mu}},& j\neq i, i\pm 1,\\
\end{cases}
\label{eq3.16x}
\end{equation}
where $i,j = 1,2,\cdots,N-1$. It is easy to see that the matrix $A$ is a real symmetric Toeplitz matrix,
which can be stored with only $(N-1)$ entries \cite{Chan96,Ng2004,Duo2018a}. Moreover, we can give the following conclusions:
\begin{proposition}
According to the definition of $A$, it holds
\begin{itemize}
\item[1)] $A$ is a strictly diagonally dominant $M$-matrix;
\item[2)] $A$ is symmetric positive definite;
\item[3)] The absolute values of the entries $a_{ij}$ away from the diagonals decay gradually, i.e.,
$a_{11} > |a_{12}| > \cdots > |a_{1,N-1}|$ and $\lim_{N\rightarrow\infty}|a_{1,N-1}| = 0$.
\end{itemize}
\label{pro3.1}
\end{proposition}
\noindent\textbf{Proof}. 1) Since $A$ is a symmetric Toeplitz matrix, then the diagonal entries are equal
to $a_{11} > 0$. Moreover, it is not hard to see that $a_{ij} < 0~(i\neq j)$. So we conclude that
$A$ is an $M$-matrix \cite[p. 533]{Horn12} and obtain
\begin{equation}
\begin{split}
a_{11} - \sum_{j\neq 1}|a_{1j}| & = a_{11} - \sum^{N-1}_{j=2}|a_{1j}| \\
& > C^{h}_{\alpha,\mu}\left[\frac{2^{\nu} + \kappa_{\mu} - 1}{2} + \frac{N^{\nu} - (N-1)^{\nu}}{N^{\mu}} + \frac{2\nu}{\alpha N^{\alpha}}\right] > 0,
\end{split}
\end{equation}
and
\begin{equation}
a_{N-1,N-1} - \sum_{j\neq N-1}|a_{N-1,j}| = a_{11} - \sum^{N-1}_{j = 2}|a_{1j}| > 0.
\end{equation}
Similarly, it follows that
\begin{equation}
\begin{split}
a_{ii} - \sum_{j\neq i}|a_{ij}| & = a_{11} - \sum_{j\neq i}a_{ij} \\
& > C^{h}_{\alpha,\mu}\left[\frac{N^{\nu} - (N-1)^{\nu}}{N^{\mu}} + \frac{2\nu}{\alpha N^{\alpha}}\right]
> 0,\quad i = 2,3,\cdots,N-2,
\end{split}
\end{equation}
a combination of the aforementioned three inequalities verifies that $A$ is a strictly diagonally dominant $M$-matrix.
2) Since $A$ is a symmetric strictly diagonally dominant $M$-matrix and all its diagonal elements are positive, i.e.,
$a_{ii} = a_{11} > 0$, then $A$ is indeed a symmetric positive definite matrix \cite[Corollary 7.2.3]{Horn12}.

3) First of all, we rewrite the matrix $A = C^{h}_{\alpha,\mu}\tilde{A} = C^{h}_{\alpha,\mu}[\tilde{a}_{ij}]_{i,j=1,\cdots,N-1
}$~--cf. Eq. (\ref{eq3.16x}. Meanwhile, it is easy to note that $a_{11} > |a_{12}|$, then we find
\begin{equation}
\begin{split}
|\tilde{a}_{12}| - |\tilde{a}_{13}| & = \frac{2^{\nu} + \kappa_{\mu} - 1}{2} - \frac{3^{\nu} - 1}{2^{\mu + 1}}\\
& \geq \frac{4^{\nu + \alpha/2} - 3^{\nu} + 1}{2^{\nu + \alpha + 1}} > 0.
\end{split}
\end{equation}
For $j = 3,4,\cdots$, we set $|j - i| = |j - 1| = k$, thus $k\geq 2$ and
\begin{equation}
\begin{split}
|\tilde{a}_{1j}| &:= f(k) = \frac{k^{-\alpha}}{2}\cdot\left[\left(1 + \frac{1}{k}\right)^{\nu} - \left(1 - \frac{1}{k}
\right)^{\nu}\right]\\
& = \left[\binom{\nu}{1}k^{-1 - \alpha} + \binom{\nu}{3}k^{-3 - \alpha} + \binom{\nu}{5}k^{-5 - \alpha} + \cdots
\right]\\
& \backsim \mathcal{O}(k^{-1 - \alpha}),
\end{split}
\end{equation}
which should implies that $f(k) > f(k + 1)$. Therefore, it follows that $a_{11} > |a_{12}| > \cdots >
|a_{1,N-1}|$. \hfill$\Box$

According to Proposition \ref{pro3.1}, if we define
\begin{equation}
\mathcal{D}(C) = \min_{1\leq i\leq N-1}\left(|C_{ii}| - \sum_{1\leq j\leq N-1,j\neq i}|C_{ij}|\right)
\end{equation}
for any matrix $C = [C_{ij}]_{i,j=1,\cdots,N-1}$, then it follows that $\mathcal{D}(A)\geq 0$, which is helpful
in the next context. The following properties of the operator $\mathcal{D}(\cdot)$ can be given as follows,
\begin{lemma}{\rm(\cite[Lemma 3]{Lin19})}
Let $C_1,C_2\in\mathbb{R}^{(N-1)\times(N-1)}$. Suppose both $C_1$ and $C_2$ have positive diagonal entries, then
it follows that $\mathcal{D}(C_1 + C_2)\geq \mathcal{D}(C_1) + \mathcal{D}(C_2)$.
\label{lem3.1}
\end{lemma}
\begin{lemma}{\rm(\cite[Lemma 4]{Lin19})}
Let $C\in\mathbb{R}^{(N-1)\times(N-1)}$. Suppose $\mathcal{D}(C)\geq 0$. Then for any nonnegative diagonal matrix
$K\in\mathbb{R}^{(N-1)\times(N-1)}$, it holds $\mathcal{D}(KC)\geq \mathcal{D}(C)\min\limits_{1\leq j\leq N-1}K_{jj}
\geq 0$.
\label{lem3.2}
\end{lemma}

Next, we exploit the above two lemmas to give the following estimation about the coefficient matrices $\mathcal{M}^{(m)}$
of Eq. (\ref{eq3.16gu}).
\begin{theorem}
For any $\frac{b^{(m,\gamma)}_m}{\Gamma(1 - \gamma)} > 0$ and $1\leq m\leq M$, it holds $\min_{1\leq m\leq M}\mathcal{D}
\Big(\frac{1}{\Gamma(1 - \gamma)}b^{(m,\gamma)}_m I +\linebreak K^{(m)}A\Big) \geq \frac{b^{(m,\gamma)}_m}{\Gamma(1 - \gamma)}$.
\label{thm3.1}
\end{theorem}
\noindent\textbf{Proof}. From Lemmas \ref{lem3.1}--\ref{lem3.2} and Proposition \ref{pro3.1}, we obtain
\begin{equation}
\begin{split}
\mathcal{D}\left(\frac{1}{\Gamma(1 - \gamma)}b^{(m,\gamma)}_m I + K^{(m)}A\right)& \geq \mathcal{D}\left(\frac{1}{\Gamma(1
- \gamma)}b^{(m,\gamma)}_m I\right) + \mathcal{D}(K^{(m)}A)\\
& \geq \frac{b^{(m,\gamma)}_m}{\Gamma(1 - \gamma)},\quad 1\leq m\leq M,
\end{split}
\end{equation}
from which the result follows. \hfill $\Box$

Before proving the final result of this section on the unconditional stability and convergence property of the FIDS
(\ref{eq3.15}), {\color{red}we recall the} following useful lemma.
\begin{lemma} Suppose $C\in\mathbb{R}^{(N - 1)\times (N-1)}$ satisfies $\mathcal{D}(C)\geq \lambda > 0$. Then, for any
${\bm y}\in\mathbb{R}^{N-1}$, it holds $\lambda\|{\bm y}\|_{\infty}\leq \|C{\bm y}\|_{\infty}$.
\label{lem3.3}
\end{lemma}
\noindent\textbf{Proof}. Since $\mathcal{D}(C)\geq \lambda > 0$, then $\mathcal{D}\left(\frac{1}{\lambda}C\right)\geq
1$ and $\|{\bm y}\|_{\infty}\leq \|\frac{C}{\lambda}{\bm y}\|_{\infty}$ \cite[Lemma 7]{Lin19}, which proves the above result.
\hfill$\Box$
\begin{theorem}
The proposed FIDS (\ref{eq3.15}) with $\epsilon\leq c_1M^{\gamma}$ is uniquely solvable and unconditionally stable in the sense that
\begin{equation}
\|{\bm u}^k\|_{\infty}\leq \|{\bm u}^0\|_{\infty} + \Gamma(1 - \gamma)\max_{1\leq s\leq k}
\frac{\|{\bm f}^{s}\|_{\infty}}{b^{(s,\gamma)}_1},\quad k = 1,2,\cdots,M,
\label{eq3.23}
\end{equation}
where $\|{\bm f}^s\|_{\infty}\leq \max_{1\leq i\leq N-1}|f^{s}_i|$.
\label{thm3.1x}
\end{theorem}
\noindent\textbf{Proof}. To prove the unique solvability is equivalent to show the invertibility of $\mathcal{M}^{(m)}$
for each $1\leq m\leq M$. By means of Theorem \ref{thm3.1} and Lemma \ref{lem3.3}, it follows that
\begin{equation}
\|\mathcal{M}^{(m)}{\bm y}\|_{\infty} = \left\|\left[\frac{1}{\Gamma(1 - \gamma)}b^{(m,\gamma)}_m I + K^{(m)}A\right]
{\bm y}\right\|_{\infty} \geq \frac{b^{(m,\gamma)}_m}{\Gamma(1 - \gamma)}\|{\bm y}\|_{\infty},\quad \forall {\bm y}\in\mathbb{R}^{N-1},
\label{eq3.26}
\end{equation}
where $1\leq m\leq M$. Therefore, $\mathcal{M}^{(m)}:~\mathbb{R}^{N-1}\mapsto \mathbb{R}^{N - 1}$ is clearly an injection for each $1\leq m\leq M$, whose
null space is simply $\{{\bm 0}\}$. Hence, $\mathcal{M}^{(m)}$'s are nonsingular, which proves the unique solvability.

On the other hand, we apply Eq. (\ref{eq3.26}) and the monotonicity of $\{b^{(m,\gamma)}_k\}~(1\leq m\leq M)$
\cite[Lemma 2.4]{Shen18} to obtain
\begin{equation*}
b^{(m,\gamma)}_m \|{\bm u}^{m}\|_{\infty}\leq \sum\limits^{m-1}_{k=1}(b^{(m,\gamma)}_{k+1} -
b^{(m,\gamma)}_k)\|{\bm u}^{k}\|_{\infty} + b^{(m,\gamma)}_1\left[\|{\bm u}^{0}\|_{\infty} + \frac{\Gamma(1 -
\gamma)}{b^{(m,\gamma)}_1}\|{\bm f}^{m}\|_{\infty}\right],~1\leq m\leq M.
\end{equation*}
Then the inequality (\ref{eq3.23}) can be proved by the method of mathematical induction, which
is similar to the proof of \cite[Theorem 4.1]{Shen18}, we omit the details here. \hfill$\Box$

On the other hand, we replace the coefficients $\{b^{(m,\gamma)}_k\}$ with $\{a^{(m,\gamma)}_k\}$
in Theorem \ref{thm3.1}, Eqs. (\ref{eq3.16gu}) and (\ref{eq3.26}), then we can obtain the following conclusion, which is helpful to analyze
the stability and convergence of DIDS (\ref{eq3.16}).
\begin{theorem}
The proposed DIDS (\ref{eq3.16}) is uniquely solvable and unconditionally stable in the sense that
\begin{equation}
\|{\bm u}^k\|_{\infty}\leq \|{\bm u}^0\|_{\infty} + \Gamma(1 - \gamma)\max_{1\leq s\leq k}
\frac{\|{\bm f}^{s}\|_{\infty}}{a^{(s,\gamma)}_1},\quad k = 1,2,\cdots,M.
\label{eq3.23x}
\end{equation}
%
\label{thm3.2}
\vspace{-3mm}
\end{theorem}
\vspace{-2mm}
\noindent\textbf{Proof}. The proof of this theorem is similar to Theorem \ref{thm3.2}, we omit
the details here. \hfill$\Box$

From Theorems \ref{thm3.1}--\ref{thm3.2}, we can see that both FIDS (\ref{eq3.15}) and DIDS (\ref{eq3.16})
are stable to the initial value $\phi$ and the right hand term $f$. Now, we consider the convergence of
these two difference schemes.
\begin{theorem}
Let $\{U^{m}_i|0 \leq i \leq N, 0 \leq m \leq M\}$ and $\{u^{m}_i|0 \leq i \leq N, 0 \leq m \leq M\}$ be,
respectively, the solutions of the problem (\ref{eq1.1}) and the difference scheme (\ref{eq3.15}). If $
\epsilon \leq \min\{c_1M^{\gamma}, T^{-\gamma/2}\}$, then
\begin{equation}
\|{\bm e}^m\|_{\infty}\leq 2\Gamma(1 - \gamma)c_2T^{\gamma}(M^{-\min\{r\gamma,2-\gamma\}} + h^p + \epsilon),
\quad 1\leq m\leq M,
\label{eq3.17}
\end{equation}
where $e^{m}_i = U^{m}_i - u^{m}_i,~0\leq i\leq N,~0\leq m\leq M$ and $p$ would be
determined by the spatial regularity of $u(x,t)$.
\label{thm3.4}
\end{theorem}
\noindent\textbf{Proof}. Writing the system (\ref{eq3.14}) as
\begin{equation*}
\begin{cases}
\frac{1}{\Gamma(1 - \gamma)}\left[b^{(m,\gamma)}_mU^{m}_i -
\sum\limits^{m-1}_{k=1}(b^{(m,\gamma)}_{k+1} - b^{(m,\gamma)}_k)U^{k}_i - b^{(m,\gamma)}_1
U^{0}_i\right] = -\kappa(-\Delta)^{\alpha/2}_{h,\mu}U^{m}_i \\
\qquad\qquad\qquad\qquad\qquad\qquad\qquad\qquad\quad~+~f^{m}_i + R^{m}_i,\quad 1\leq i\leq N-1,~1\leq m\leq M,\\
U^{m}_0 = U^{m}_{N} = 0,\qquad\qquad\qquad\qquad\qquad\qquad\qquad\qquad~~\quad~0\leq m \leq M,\\
U^{0}_i = \phi(x_i),\qquad\qquad\qquad\qquad\qquad\qquad\qquad\qquad~~\qquad~~1\leq i\leq N - 1,
\end{cases}
\end{equation*}
and subtracting the Eq. (\ref{eq3.15}) from the corresponding above system
\begin{equation}
\begin{cases}
\frac{1}{\Gamma(1 - \gamma)}\left[b^{(m,\gamma)}_me^{m}_i -
\sum\limits^{m-1}_{k=1}(b^{(m,\gamma)}_{k+1} - b^{(m,\gamma)}_k)e^{k}_i - b^{(m,\gamma)}_1
e^{0}_i\right] = -\kappa(-\Delta)^{\alpha/2}_{h,\mu}e^{m}_i \\
\qquad\qquad\qquad\qquad\qquad\qquad\qquad\qquad\quad~ + R^{m}_i,\quad 1\leq i\leq N-1,~1\leq m\leq M,\\
e^{m}_0 = e^{m}_{N} = 0,\qquad\qquad\qquad\qquad\qquad\qquad\qquad\quad~~0\leq m \leq M,\\
e^{0}_i = 0,\qquad\qquad\qquad\qquad\qquad\qquad\qquad\qquad\qquad~1\leq i\leq N - 1,
\end{cases}
\end{equation}
By means of Theorem \ref{thm3.1} and the matrix analysis described above, it follows that
\begin{equation*}
\|{\bm e}^m\|_{\infty}\leq \Gamma(1 - \gamma)\max_{1\leq s\leq m}
\frac{\|{\bm R}^{s}\|_{\infty}}{b^{(s,\gamma)}_1},\quad m = 1,2,\cdots,M,
\end{equation*}
the rest of this proof is also similar to \cite[Theorem 4.2]{Shen18}.

Again, we employ the similar strategy to give the error analysis of DIDS (\ref{eq3.16}) as follows.
\begin{theorem}
Let $\{U^{m}_i|0 \leq i \leq N, 0 \leq m \leq M\}$ and $\{u^{m}_i|0 \leq i \leq N, 0 \leq m \leq M\}$ be,
respectively, the solutions of the problem (\ref{eq1.1}) and the difference scheme (\ref{eq3.16}), then
\begin{equation}
\|{\bm e}^m\|_{\infty}\leq \mathcal{O}(M^{-\min\{r\gamma,2-\gamma\}} + h^p),
\quad 1\leq m\leq M,
\end{equation}
where $e^{m}_i = U^{m}_i - u^{m}_i,~0\leq i\leq N,~0\leq m\leq M$.
\label{thm3.5}
\end{theorem}
\vspace{-2mm}
In practice, the value of $\epsilon$ is sufficiently small such that the tolerance
error in (\ref{eq3.17}) can be negligible in compared to the space and time errors. Then it also finds that the
numerical errors for DIDS and FIDS are almost identical but the later is often faster -- cf. Section \ref{sec5}. With the help of arguments in proving \cite[Theorems 3.1-3.2]{Duo2018a}, it is not hard to make the above
convergence results described in Theorems \ref{thm3.4}--\ref{thm3.5} more specific.
\begin{remark} For determining the value of $p$, it reads
\begin{itemize}
\item Suppose that the solution of the problem (\ref{eq1.1}) satisfies the condition (\ref{eq2.1}) and $u(x,
\cdot) \in \mathcal{C}^{s,\alpha/2}(\mathbb{R})$ with $(x,t) \in \mathbb{R}\times[0,T]$ and $s \geq 1$, then
solutions of FIDS (\ref{eq3.15}) and DIDS (\ref{eq3.16}) with $\mu\in(\alpha,2)$ converge to the exact solutions
of Eq. (\ref{eq1.1}), respectively;
\item For $s = 1$, $\alpha\in(0,2)$ and $\mu\in(\alpha,2]$, the convergence rates of FIDS (\ref{eq3.15}) and
DIDS (\ref{eq3.16}) are (at least) $\mathcal{O}(M^{-\min\{r\gamma,2-\gamma\}} + h^{1 - \frac{\alpha}{2}}+ \epsilon)$ and $\mathcal{
O}(M^{-\min\{r\gamma,2-\gamma\}} + h^{1 - \frac{\alpha}{2}})$, respectively; see \cite{Duo2018a,Duo2019z} for details.
\item For $s \geq 3$ and $\alpha\in(0,2)$, the convergence rates of FIDS (\ref{eq3.15}) and DIDS (\ref{eq3.16})
with $\mu = 2$ or $\mu = 1 + \alpha/2$ are $\mathcal{O}(M^{-\min\{r\gamma,2-\gamma\}} + h^2 + \epsilon)$
and $\mathcal{O}(M^{-\min\{r\gamma,2-\gamma\}} + h^2)$, respectively.
\end{itemize}
\label{rem3.1}
\end{remark}
Moreover, we will work out some numerical results for supporting the above theoretical
convergence behaviors described in Section \ref{sec5}.
\section{Efficient implementation based on preconditioning of the difference schemes}
\label{sec4}
In the section, we analyze both the implementation and computational complexity of FIDS (\ref{eq3.15})
and DIDS (\ref{eq3.16}) and we propose an efficient implementation utilized preconditioned Krylov subspace
solvers. Noting that $a^{ (m, \gamma) }_m = b^{ (m,\gamma) }_m > 0$, we start the efficient implementation
from the following matrix form of these two implicit difference schemes at the time level $1\leq m\leq M$,
which are given by Eq. (\ref{eq3.16gu}) and
\vspace{-2mm}
\begin{equation}
\mathcal{M}^{(m)}{\bm u}^{m} = \frac{1}{\Gamma(1 - \gamma)} \left[\sum\limits^{m-1}_{k=1}(a^{(m,\gamma)}_{k+1}
- a^{(m,\gamma)}_k){\bm u}^{k} + a^{(m,\gamma)}_1{\bm u}^{0}\right] + {\bm f}^{m},
\label{eq4.1}
\end{equation}
respectively. From Theorems \ref{thm3.1x}-\ref{thm3.2}, it knows that both Eq. (\ref{eq3.16gu}) and Eq.
(\ref{eq4.1}) have the unique solutions. In addition, it is meaningful to remark that Eq. (\ref{eq3.16gu}) and Eq.
(\ref{eq4.1}) corresponding to FIDS (\ref{eq3.15}) and DIDS (\ref{eq3.16}) are inherently sequential, thus they
are both difficult to parallelize over time.
\vspace{-2mm}
\subsection{The circulant preconditioner}
On the other hand, it is useful to note that the matrix-vector product $\mathcal{M}^{(m)}{\bm v}$ can be efficiently
calculated by
\vspace{-2mm}
\begin{equation}
\mathcal{M}^{(m)}{\bm v} = \frac{1}{\Gamma(1 - \gamma)}b^{(m,\gamma)}_m{\bm v} + K^{(m)}(A{\bm v}),
\end{equation}
where ${\bm v}\in\mathbb{R}^{N-1}$ is any vector and the Toeplitz matrix-vector $A{\bm v}$ can be implicitly
evaluated via the FFTs in $\mathcal{O}(N\log N)$ operations. In other words, we can
use a \textit{matrix-free} method to compute $\mathcal{M}^{(m)}{\bm v}$ quickly. Based on such observations,
the Krylov subspace method should be the most suitable solver for Eq. \ref{eq3.16gu} or Eq. (\ref{eq4.1}) one
by one. However, when the coefficients and the order of integral fractional Laplacian are not small, then the
coefficient matrices $\mathcal{M}^{(m)} $ will be increasingly ill-conditioned (cf. Section \ref{sec5}). This fact deeply slows down the convergence of the Krylov subspace method, while the preconditioning techniques are often used to overcome this
difficulty \cite{Lei2013,Gu2017,Chan96,Ng2004}. In the literature on Toeplitz systems, circulant preconditioners always played important roles
\cite{Chan96,Ng2004}. In fact, circulant preconditioners have been theoretically and numerically studied with
applications to fractional partial differential equations for recent years; see for instance \cite{Lei2013,Gu2017,Hao2019}.

In this work, we design a family of the Strang's preconditioners \cite{Chan96} for accelerating the convergence of Krylov subspace
solvers. More precisely, the circulant preconditioners are given for Eq. \ref{eq3.16gu} and/or Eq. (\ref{eq4.1})
as follows,
\begin{equation}
\mathcal{P}^{(m)} =
\frac{1}{\Gamma(1 - \gamma)}b^{(m,\gamma)}_mI + \kappa^{(m)} s(A) = F^{*}\left[\frac{1}{\Gamma(1
- \gamma)}b^{(m,\gamma)}_m + \kappa^{(m)}\Lambda\right]F,
\label{eq4.3}
\end{equation}
where $F$ and $F^{*}$ are the Fourier matrix and its conjugate transpose, respectively, and the scalar $\kappa^{(m)}
= \frac{1}{N-1}\sum^{N-1}_{i = 1}\kappa^{m}_i$. Meanwhile, $s(A) = F^{*}\Lambda F$ is the Strang circulant
approximation \cite{Chan96,Ng2004} of the Toeplitz matrix $A$ and the diagonal matrix $\Lambda$ contains all the
eigenvalues of $s(A)$ with the first column: ${\bm c}_{S} = [a_{11}, \cdots, a_{1,\lfloor\frac{N+1}{2}\rfloor},
a_{1,\lfloor\frac{N}{2}\rfloor},\cdots,a_{12}]^{\top}\in\mathbb{R}^{N-1}$.
Therefore the matrix $\Lambda = {\rm diag}(F{\bm c}_{S})$ can be {\em computed in advance and only one time during each
time level}. Besides, as $\mathcal{P}^{(m)}$ are the circulant matrices, we observe from Eq. (\ref{eq4.3}) that the
inverse-matrix-vector product ${\bm z} = [\mathcal{P}^{(m)}]^{-1}{\bm v}$ can be carried out in $\mathcal{O}(N\log
N)$ operations via the (inverse) FFTs. In one word, we exploit a fast preconditioned Krylov subspace method with only $\mathcal{O}(N)$ memory requirement and $\mathcal{O}(N\log
N)$ computational cost per iteration, while the number of iterations and the computational cost are greatly
reduced.

To investigate the properties of the proposed preconditioners, the following lemma is the key to prove the invertibility
of $\mathcal{P}^{(m)}$ in Eq. (\ref{eq4.3}).
\vspace{-2mm}
\begin{lemma}
All eigenvalues of $s(A)$ fall inside the open disc
\begin{equation}
\{z\in\mathbb{R}:~|z - a_{11}| < a_{11}\},
\end{equation}
and all the eigenvalues of $s(A)$ are strictly positive for all $N$.
\label{lem4.1}
\end{lemma}
\noindent\textbf{Proof}. First of all, since the matrix $A$ is symmetric, then $s(A)$ is also symmetric and its eigenvalues
should be real. All the Gershgorin disc of the circulant matrix $s(A)$ are centered at $a_{11}$ with radius
\begin{equation}
r_N = 2\sum^{\lfloor\frac{N+1}{2}\rfloor}_{\ell = 2}|a_{1,\ell}| < a_{11},
\end{equation}
the above inequality holds due to the expression of Eq. (\ref{eq3.16x}), where the expression of $a_{11}$ contains exactly
the sum of $2|a_{1,\ell}|~(\ell = 2,3,\ldots,N-1)$. In conclusion, all the eigenvalues of $s(A)$ are
strictly positive for all $N$. \hfill $\Box$

According to Lemma \ref{lem4.1}, it means that $s(A)$ is a real symmetric positive definite matrix. Moreover, the invertibility
of circulant preconditioners $\mathcal{P}^{(m)}$ (\ref{eq4.3}) can be given for all $m = 1,2,\cdots,M$ as follows.
\begin{lemma}
Let $\alpha\in(0,2)$. The preconditioner $\mathcal{P}^{(m)}$ is invertible and
\begin{equation}
\left\|\left(\mathcal{P}^{(m)}\right)^{-1}\right\|_2 <
\frac{1}{\frac{b^{(m,\gamma)}_m}{\Gamma(1 -\gamma)}~+~\kappa^{(m)}\cdot\min\limits_{1\leq
k\leq N-1}[\Lambda]_{k,k}},
\vspace{-3mm}
\end{equation}
\label{lem4.2}
\end{lemma}
\noindent\textbf{Proof}. According to Lemma \ref{lem4.1}, we have $[\Lambda]_{k,k} > 0$. Noting that $a^{(m,\gamma)}_m > 0$
(or $b^{(m,\gamma)}_m > 0$) and $\kappa^{(m)} > 0$, we have
\begin{equation}
\left[\frac{1}{\Gamma(1 - \gamma)}b^{(m,\gamma)}_m + \kappa^{(m)}\Lambda\right]_{k,k} > 0
\end{equation}
for $k = 1,2,\cdots,N-1$. Therefore, $\mathcal{P}^{(m)}$ is invertible. Furthermore, we have
\begin{equation*}
\left\|\left(\mathcal{P}^{(m)}\right)^{-1}\right\|_2 = \frac{1}{\min\limits_{1\leq k\leq N-1}\left[\frac{1}{\Gamma(1 -
\gamma)}b^{(m,\gamma)}_m + \kappa^{(m)}\Lambda\right]_{k,k}} < \frac{1}{\frac{b^{(m,\gamma)}_m}{\Gamma(1 -\gamma)}
+ \kappa^{(m)}\cdot\min\limits_{1\leq k\leq N-1}[\Lambda]_{k,k}},
\end{equation*}
the other inequality can be similarly obtained. \hfill $\Box$
\subsection{Spectrum of the preconditioned matrix}
\label{sec4.2}
In this subsection, we study the spectrum of the preconditioned matrix, which can help us to understand the
convergence of preconditioned Krylov subspace solvers. For convenience of our investigation, we first
assume that the diffusion coefficient function $\kappa(x,t) \equiv \kappa(t)$, then Eq.~(\ref{eq3.16gu}) or
Eq.~(\ref{eq4.1}) will be a sequence of real symmetric positive definite linear systems, where the coefficient
matrices reduce to $\mathcal{M}^{(m)} = \frac{1}{\Gamma(1 - \gamma)}b^{(m,\gamma)}_mI + \kappa^{(m)}A$ corresponding
to Eq. (\ref{eq3.16gu}) and Eq.~(\ref{eq4.1}), respectively. The preconditioned CG (PCG) method \cite{Ng2004} should be a suitable
candidate for solving such linear systems one by one. Moreover, the spectrum of the preconditioned matrix
$\left(\mathcal{P}^{(m)}\right)^{-1}\mathcal{M}^{(m)}$ are available for both Eq. (\ref{eq3.16gu}) and Eq.~(\ref{eq4.1})
at each time level $m$, so we take Eq. (\ref{eq3.16gu}) as the research object in the next context.

Throughout this subsection, we rewrite Eq. (\ref{eq3.16gu}) into the following equivalent form
\begin{equation}
\tilde{\mathcal{M}}^{(m)}{\bm u}^{m} = \frac{1}{\Gamma(1 - \gamma)C^{h}_{\alpha,\mu}} \left[\sum\limits^{m-1}_{k=1}
(b^{(m,\gamma)}_{k+1} - b^{(m,\gamma)}_k){\bm u}^{k}~{\color{red}+}~ b^{(m,\gamma)}_1{\bm u}^{0}\right]
+ \frac{1}{C^{h}_{\alpha,\mu}}{\bm f}^{m},
\label{eq4.9}
\end{equation}
where $\tilde{\mathcal{M}}^{(m)} = \frac{b^{(m,\gamma)}_m}{\Gamma(1 - \gamma)C^{h}_{\alpha,\mu}}I + \kappa^{(m)}\tilde{A}$
and its corresponding circulant preconditioner $\mathcal{P}^{(m)}$ reduces to $\tilde{\mathcal{P}}^{(m)} = \frac{b^{(m,\gamma)}_m}{\Gamma(1
- \gamma)C^{h}_{\alpha,\mu}}I + \kappa^{(m)}s(\tilde{A})$, which is still invertible -- cf. Lemma \ref{lem4.2}. Moreover, we assume that $M$ and $r$ are properly chosen, depending on
$N$, such that $\eta^{(m)}_{N,M,r} \triangleq \frac{b^{(m,\gamma)}_m}{\Gamma(1 - \gamma)C^{h}_{\alpha,\mu}}$ in (\ref{eq4.9}) is
bounded away from 0; i.e., there exist two real numbers ``$\check{\eta}$" and ``$\hat{\eta}$" such that
\begin{equation}
0 < \check{\eta}\leq \eta^{(m)}_{N,M,r}\leq \hat{\eta},\quad \forall N~{\rm and}~m = 1,2,\cdots,M-1.
\label{eq4.10}
\end{equation}
We add a subscript $N$ to each matrix to denote the matrix size. Under the above assumption in (\ref{eq4.10}), the matrix
$\tilde{\mathcal{M}}^{(m)}$, $K^{(m)}$, $\eta^{(m)}_{N,M,r}$, and $\tilde{\mathcal{P}}^{(m)}$ are independent of $m$, and
we therefore simply denote them as $\tilde{\mathcal{M}}_{N-1}$, $\kappa$ (constant), $\eta_{N,M,r}$ (constant) and $\tilde{
\mathcal{P}}_{N-1}$, respectively. Now the coefficient matrix $\tilde{\mathcal{M}}^{(m)}$ in (\ref{eq3.16gu})
becomes
\begin{equation}
\begin{split}
\mathcal{M}_{N-1} & = \eta_{N,M,r} I + \kappa \tilde{A}_{N-1}
\\
& = [\eta_{N,M,r} + \kappa(\tilde{a}_{11}-\varrho)]I + \begin{bmatrix}
\varrho    &\kappa \tilde{a}_{12}    &\cdots & \kappa \tilde{a}_{1,N-2} & \kappa \tilde{a}_{1,N-1} \\
\kappa \tilde{a}_{12}    &\varrho    &\kappa \tilde{a}_{12} & \cdots    & \kappa \tilde{a}_{1,N-2} \\
\vdots    &\kappa \tilde{a}_{12}    &\varrho & \ddots    &\vdots     \\
\kappa \tilde{a}_{1,N-2} &\cdots    &\ddots & \ddots    &\kappa \tilde{a}_{12}    \\
\kappa \tilde{a}_{1,N-1} &\kappa\tilde{a}_{1,N-2} &\cdots &\kappa \tilde{a}_{12}    & \varrho
\end{bmatrix}\\
& \triangleq [\eta_{N,M,r} + \kappa(\tilde{a}_{11}-\varrho)]I + G_{N-1},
\end{split}
\label{eq4.9x}
\end{equation}
where we set $|\tilde{a}_{12}| < \varrho < \tilde{a}_{11}$ (without loss of generality), $G_{N-1} = [\tilde{g}_{ij}]_{(N-1)\times(N-1)} = [\tilde{g}_{|i-j|}]_{(N-1)\times(N-1)}$ and $\tilde{g}_{ij} = \kappa \tilde{a}_{ij}$.

To study the spectrum of the preconditioned matrix $\left(\mathcal{P}^{(m)}\right)^{-1}\mathcal{M}^{(m)}$, we first introduce
the \textit{generating function} of the sequence of Toeplitz matrices $\{G_{N}\}^{\infty}_{N=1}$:
\begin{equation}
p(\theta) = \sum^{\infty}_{k = -\infty}\tilde{g}_k e^{\iota k\theta},
\end{equation}
where $\tilde{g}_k$ is the $k$-th diagonal of $G_{N} = [\tilde{g}_{i-j}]_{N\times N}$ and $\iota = \sqrt{-1}$. The generating function $p(\theta)
$ is in the \textit{Wiener class} \cite{Chan96,Ng2004} if and only if
\begin{equation}
\sum^{\infty}_{k = -\infty}|\tilde{g}_k| < \infty.
\end{equation}
For $G_{N-1}$ defined in (\ref{eq4.9}), we have the following conclusion.
\begin{lemma}
Under the above assumptions, it finds that $p(\theta)$ is real-valued and in the Wiener class.
\label{lem4.3}
\end{lemma}
\noindent\textbf{Proof}. For convenience of our investigation, we can rewrite $p(\theta)$ for the matrix $G_{N}$
defined in (\ref{eq4.9}) as
\begin{equation}
p(\theta) = \tilde{g}_0 + 2\sum^{\infty}_{k=1}\tilde{g}_{k}\cos(k\theta) = \varrho - 2\sum^{\infty}_{k=1}(-\tilde{g}_{k})\cos(k\theta).
\end{equation}
Since it knows that $0 < -\tilde{g}_{k+1} < -\tilde{g}_{k}$, $\lim_{k\rightarrow\infty}(-\tilde{g}_k) = 0$ and
\begin{equation*}
|(-\tilde{g}_k)\cos(k\theta)| < (-\tilde{g}_k)
\end{equation*}
with the series $\sum^{\infty}_{k=1}(-\tilde{g}_k)$ being convergent, thus the series $\sum^{\infty}_{k=1}(-\tilde{g}_{k})
\cos(k\theta)$ converges to a real-valued function for $\forall\theta\in[-\pi,\pi]$, which also implies that $p(\theta)$
is real-valued. According to Proposition \ref{pro3.1} and its proof, it is not hard to note that $\lim_{k\rightarrow\infty}
|\tilde{g}_{k}| = 0$. Therefore, it follows that $\sum^{\infty}_{k=-\infty}|\tilde{g}_{k}| < \infty$, which completes the  proof.
\hfill $\Box$

In fact, Lemma \ref{lem4.3} ensures the following property that the given Toeplitz matrix $G_N$ can be approximated via a circulant matrix
well.
\begin{lemma}
If $p(\theta)$, the generating function of $G_N$, is in the Wiener class,
then for any $\epsilon > 0$, there exist $N'$ and $M' > 0$, such that for all $N > N'$,
\begin{equation}
G_N - s(G_N) = U_N + V_N,
\end{equation}
where ${\rm rank}(U_N)\leq M'$ and $\|V_N\|_2 < \epsilon$.
\label{lem4.4}
\end{lemma}

Now we consider the spectrum of $(\tilde{\mathcal{P}}_{N-1})^{-1}\tilde{\mathcal{M}}_{N-1} - I$ is clustered around 1.
\begin{theorem}
If $\eta_{N,M,r}$ satisfies the assumption (\ref{eq4.10}), for any $0< \epsilon <1$, there exists $N'$
and $M'>0$ such that, for all $N > N'$, at most $2M'$ eigenvalues of the matrix $\tilde{\mathcal{M}}_{
N-1} -\tilde{\mathcal{P}}_{N-1}$ have absolute values exceeding $\epsilon$.
\end{theorem}
\noindent\textbf{Proof}. With the help of Lemma \ref{lem4.4}, we note that
\begin{equation}
\begin{split}
\tilde{\mathcal{M}}_{N-1} -\tilde{\mathcal{P}}_{N-1} & = \kappa \tilde{A} - \kappa s(\tilde{A})\\
& = G_{N-1} - s(G_{N-1})\\
& = U_{N-1} + V_{N-1}.
\end{split}
\end{equation}
Since both $V_{N-1}$ and $U_{N-1}$ are real symmetric with $\|V_{N-1}\|_2 < \epsilon$ and ${\rm rank}(U_{
N-1}) < M'$, hence the spectrum of $V_{N - 1}$ lies in $(-\epsilon,\epsilon)$. By the celebrated Weyl's
theorem, we see that at most $2M'$ eigenvalues of $\tilde{\mathcal{M}}_{N-1} -\tilde{\mathcal{P}}_{N-1}$
have absolute values exceeding $\epsilon$. \hfill $\Box$

At this stage, we can see from Lemma \ref{lem4.2} that
\begin{equation}
\|(\tilde{\mathcal{P}}_{N-1})^{-1}\|_2 = \frac{1}{\min\limits_{1\leq k\leq N-1}\left[\eta_{N,M,r} + \tilde{
\Lambda}\right]_{k,k}} < \frac{1}{\eta_{N,M,r}}\leq \frac{1}{\check{\eta}},
\end{equation}
where $s(\kappa\tilde{A})= F^{*}\tilde{\Lambda}F$ and the diagonal matrix $\tilde{\Lambda}$ contains all
the eigenvalues of $s(\kappa\tilde{A})$. Meanwhile, we employ the fact that
\begin{equation}
(\tilde{\mathcal{P}}_{N-1})^{-1}\tilde{\mathcal{M}}_{N-1} - I  = (\tilde{\mathcal{P}}_{N-1})^{-1}U_{N-1}
- (\tilde{\mathcal{P}}_{N-1})^{-1}V_{N-1},
\end{equation}
then we have the following corollary.
\begin{corollary}
If $\eta_{N,M,r}$ satisfies the assumption (\ref{eq4.10}), for any $0< \epsilon <1$, there exists $N'$
and $M'>0$ such that, for all $N > N'$, at most $2M'$ eigenvalues of the matrix $(\tilde{\mathcal{P}}_{
N-1})^{-1}\tilde{\mathcal{M}}_{N-1} - I$ have absolute values exceeding $\epsilon$.
\end{corollary}

Thus the spectrum of $(\tilde{\mathcal{P}}_{N-1})^{-1}\tilde{\mathcal{M}}_{N-1}$ is clustered around 1 for
enough large $N$. It follows that the convergence rate of the PCG method is superlinear; refer to \cite{Chan96,Ng2004}
for details. Based on such observations, the preconditioner $\mathcal{P}^{(m)}$ is fairly predictable to
accelerating the convergence of PCG for solving both Eq. (\ref{eq3.16gu}) and Eq.~(\ref{eq4.1}) at each time level $m
=1,2,\cdots,M$ well, respectively; refer to numerical results in the next section.

Besides, although the theoretical analysis in Section \ref{sec4.2} is only available for handling the
model problem (\ref{eq1.1}) with time-varying diffusion coefficients, i.e., $\kappa(x,t) \equiv \kappa(t)$,
the preconditioner $\mathcal{P}^{(m)}$ is still efficient to accelerate the convergence of nonsymmetric
Krylov subspace solvers for Eq. (\ref{eq3.16gu}) and/or Eq.~(\ref{eq4.1}) corresponding to the
problem (\ref{eq1.1}). The variable diffusion coefficients and nonsymmetric discretized
linear systems make it greatly challenging to theoretically study the eigenvalue distributions of preconditioned
matrices $(\mathcal{P}^{(m)})^{-1}\mathcal{M}^{(m)}$, but we provide numerical results to show the clustering
eigenvalue distributions of some specified preconditioned matrices in Section \ref{sec5}. In summary, we can
analyze the computational complexity and memory requirement for both FIDS and DIDS as follows.
\begin{proposition}
The FIDS (or DIDS) has $\mathcal{O}(NN_{exp})$ (or $\mathcal{O}(NM)$) memory requirement and $\mathcal{O}(MN(
\log N + N_{exp}))$ (or $\mathcal{O}(MN(\log N + M))$) computational complexity.
\end{proposition}
%
%
%
\section{Numerical experiments}
\label{sec5}
In this section, numerical experiments are presented to achieve our two-fold
objective\footnote{Here we note that the MATLAB codes of all the numerical tests are available from
authors' emails and will make them public in the GitHub repository: \url{https://github.com/Hsien-Ming-Ku/Group-of-FDEs}.}. We show that the
proposed FIDS and DIDS can indeed converge with the theorectical accuracy in both space and time. Meanwhile,
we assess the computational efficiency and theorectical results on circulant preconditioners described in
Section \ref{sec4}. For the Krylov subspace method and direct solver, we exploit built-in functions for the
preconditioned BiCGSTAB (PBiCGSTAB) method \cite{Vorst} (in Example 2) and MATLAB's backslash in Examples 1--2,
respectively. For the BiCGSTAB method with circulant preconditioners, the stopping criterion of those methods is
$\|{\bm r}^{(k)}\|_2/\|{\bm r}^{(0)}\|_2 < {\rm tol} = 10^{-10}$, where ${\bm r}^{(k)}$ is the residual vector
of the linear system after $k$ iterations, and the initial guess is chosen as the zero vector. All experiments
were performed on a Windows 10 (64 bit) PC-Intel(R) Core(TM) i5-8265U CPU (1.6 $\sim$ 3.9 GHz), 8 GB of RAM using MATLAB 2017b
with machine epsilon $10^{-16}$ in double precision floating point arithmetic. The computing time reported is an average over 20 runs of our algorithms. We also choose the tolerance error $\epsilon = 10^{-10},~10^{-9}$ for FIDS
in Examples 1--2, respectively. Moreover, some notations on numerical errors are introduced as follows:
\[
{\rm Error}_{\infty}(N,M) = \max_{0\leq j\leq M}\|{\bm E}^{j}\|_{\infty}\quad {\rm and}\quad {\rm Error}_2(N,M) =
\max_{0\leq j\leq M}\|{\bm E}^{j}\|_2,
\]
then
\begin{equation*}
{\rm Rate}_{\infty} =
\begin{cases}
\log_2\left(\frac{{\rm Error}_{\infty}(N(M/2),M/2)}{{\rm Error}_{\infty}(N(M),M)}\right),&{\rm(temporal~convergence~order)},\\
\log_2\left(\frac{{\rm Error}_{\infty}(N/2,M(N/2))}{{\rm Error}_{\infty}(N,M(N))}\right),&{\rm(spatial~convergence~order)},
\end{cases}
\end{equation*}
and
\begin{equation*}
{\rm Rate}_2 =
\begin{cases}
\log_2\left(\frac{{\rm Error}_2(N(M/2),M/2)}{{\rm Error}_2(N(M),M)}\right),&{\rm(temporal~convergence~order)},\\
\log_2\left(\frac{{\rm Error}_2(N/2,M(N/2))}{{\rm Error}_2(N,M(N))}\right),&{\rm(spatial~convergence~order)}.
\end{cases}
\end{equation*}

\begin{figure}
\centering
\includegraphics[width=3.15in,height=2.95in]{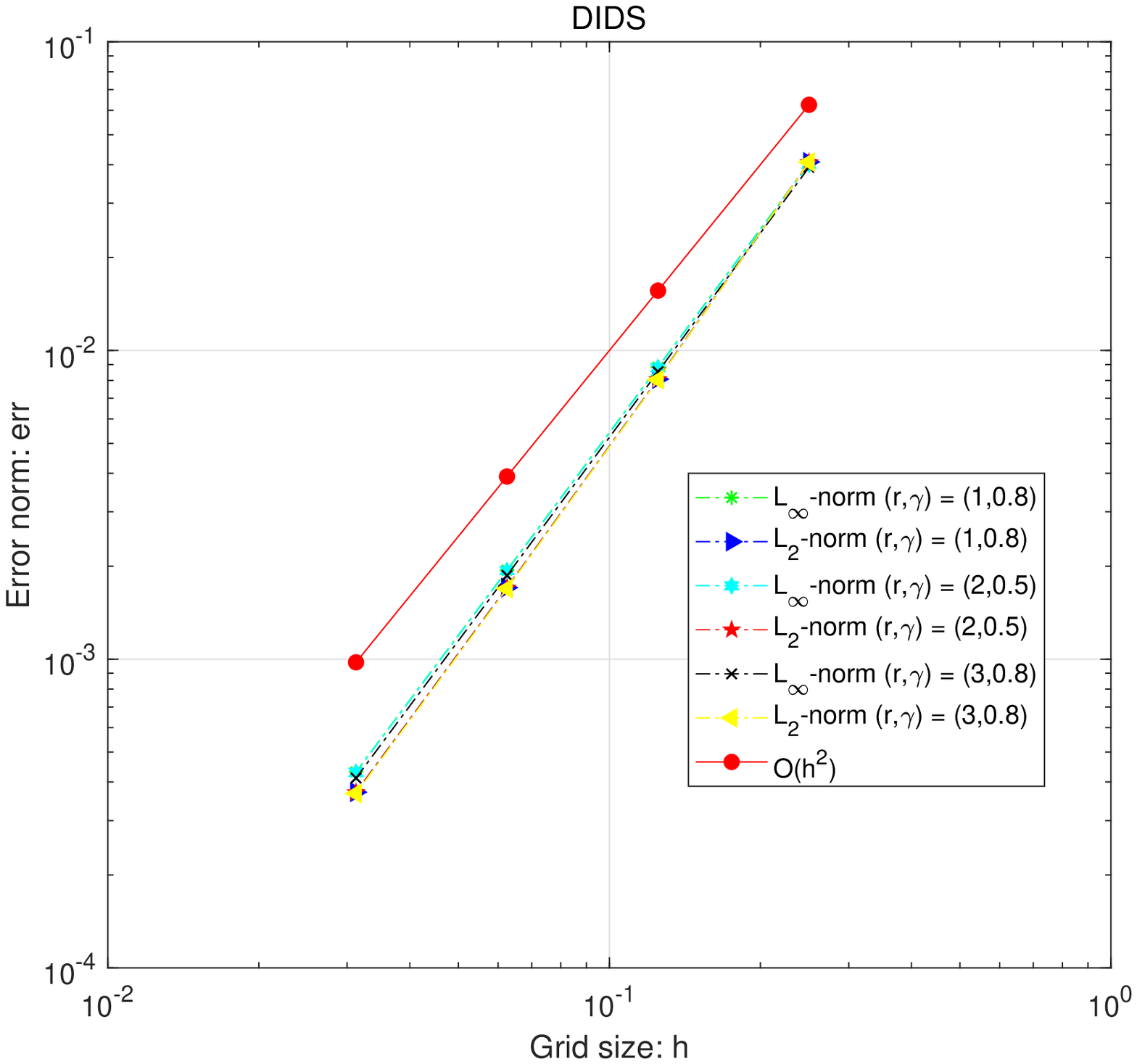}
\includegraphics[width=3.12in,height=2.95in]{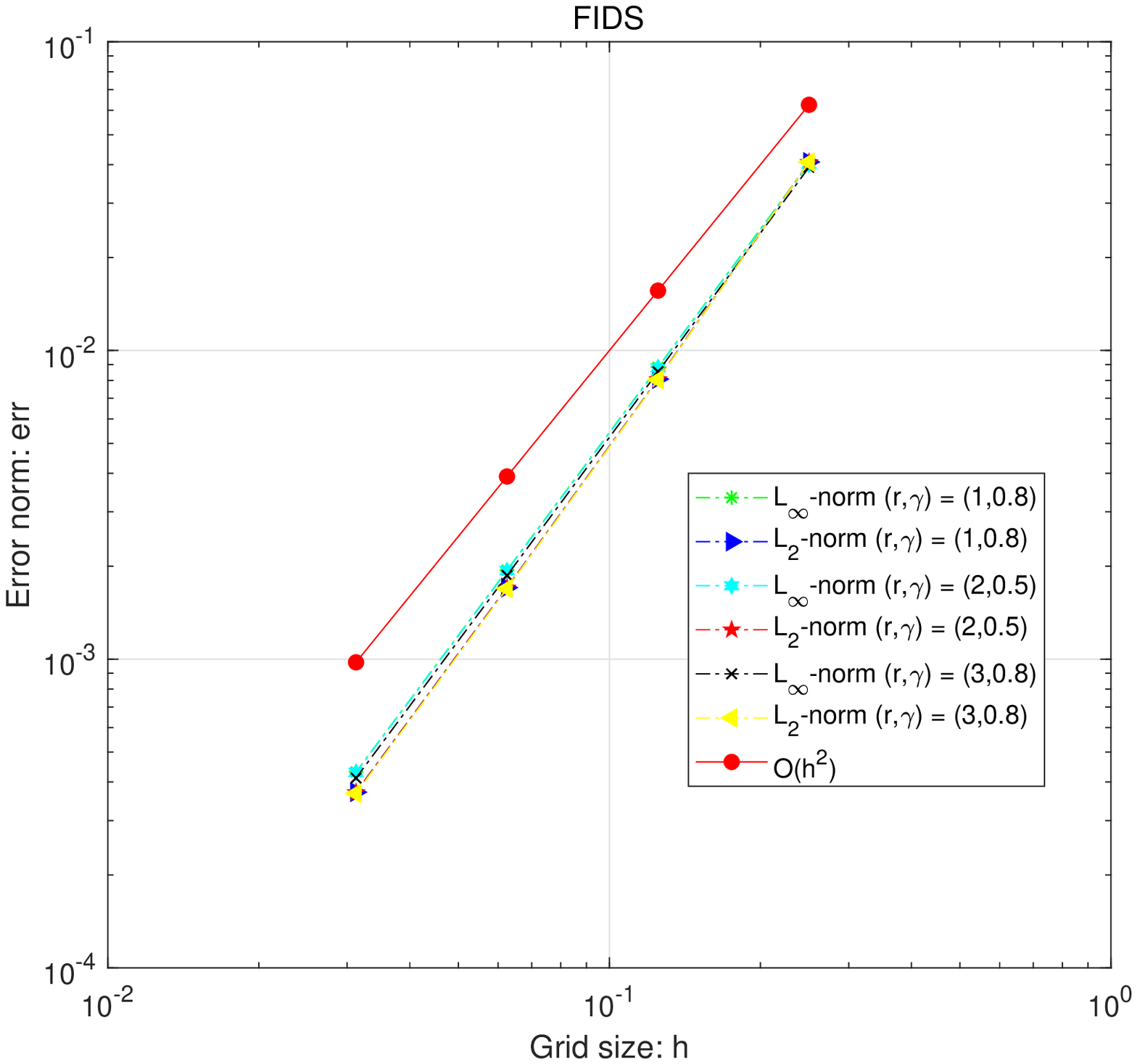}
\caption{The spatial convergence order of two difference schemes for Example 1 with $\alpha = 1.6$, $M(N) = \lceil (N/2)^{2/\min\{r\gamma,2-\gamma\}}
\rceil$, $\mu = 2$, $\kappa_{\mu} = 2$, and $s = 3$. Left: DIDS (\ref{eq3.16});~~Right: FIDS (\ref{eq3.15}).}
\label{fig1}
\end{figure}

\begin{figure}
\centering
\includegraphics[width=3.15in,height=2.95in]{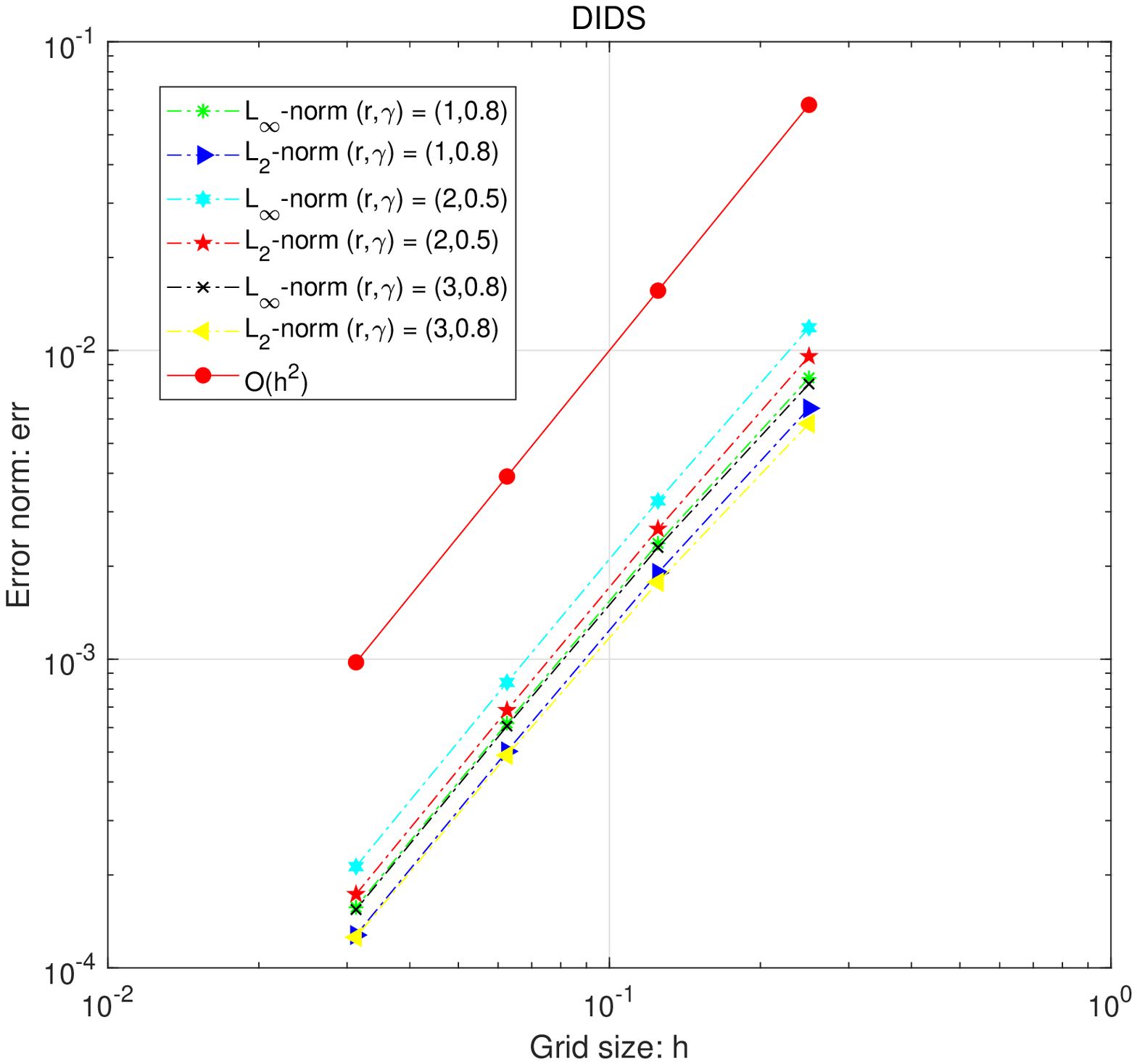}
\includegraphics[width=3.15in,height=2.95in]{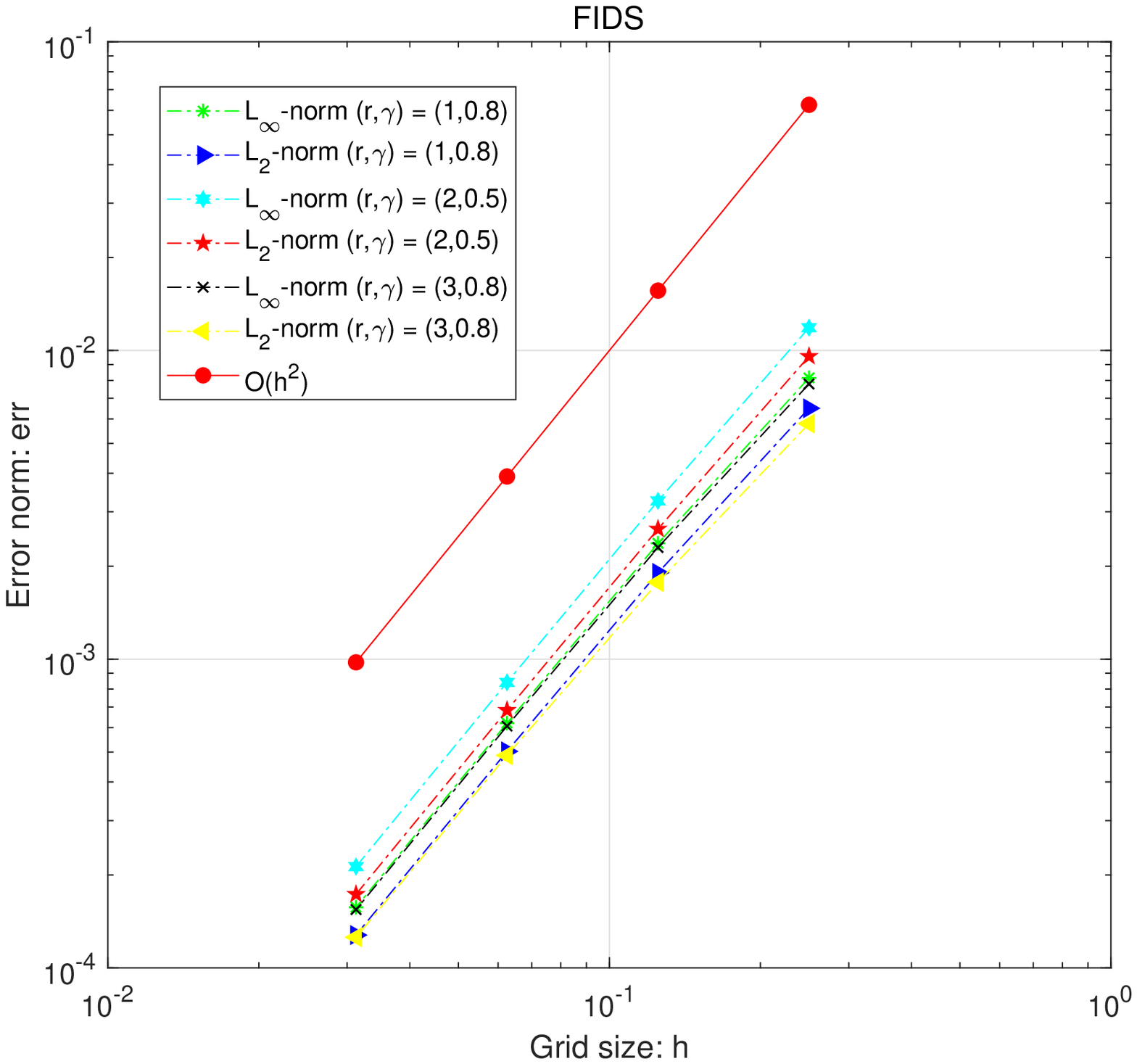}
\caption{The spatial convergence order of two difference schemes for Example 1 with $\alpha = 0.4$, $M(N) = \lceil (N/2)^{2/\min\{r\gamma,2-\gamma\}}
\rceil$, $\mu = 2$, $\kappa_{\mu} = 2$, and $s = 3$. Left: DIDS (\ref{eq3.16});~~Right: FIDS (\ref{eq3.15}).}
\label{fig2}
\end{figure}

\noindent\textbf{Example 1}. (Accuracy test) In this example, we consider Eq. (\ref{eq1.1}) with the spatial domain $\Omega =
(-1,1)$ and the time interval $[0,T] = [0,1]$. The diffusion coefficients $\kappa(x,t) = (1 + t)e^{0.8x + 1}$ and
the source term is given
\begin{equation*}
\begin{split}
f(x,t) & = \Gamma(1 + \gamma)(1 - x^2)^{s + \alpha/2} + \kappa(x,t)
\frac{2^{\alpha}\Gamma(\frac{\alpha + 1}{2})\Gamma(s + 1 + \alpha/2)}{\sqrt{\pi}\Gamma(s + 1)}~\times\\
&~~~{}_2F_1\left(\frac{\alpha + 1}{2},-s;\frac{1}{2};x^2\right)(t^{\gamma} + 1),\quad s\in\mathbb{N}^{+},
\end{split}
\end{equation*}
where ${}_2F_1(\cdot)$ denotes the Gauss hypergeometric function which can be computed via the MATLAB built-in
function `{\tt hypergeom.m}' and the initial-boundary value conditions are
\begin{equation*}
u(x,t) = 0,~x\notin\Omega,\quad {\rm and}\quad u(x,0) = (1 - x^2)^{s + \alpha/2}.
\end{equation*}
Thus the exact solution of this problem is $u(x,t) = (1 - x^2)^{s + \alpha/2}(t^{\gamma} + 1)$. The numerical
results involving both spatial and temporal convergence orders of FIDS (\ref{eq3.15}) and DIDS (\ref{eq3.16})
are shown in Tables \ref{tab1}--\ref{tab6} and Figs. \ref{fig1}--\ref{fig2}. Here it should mentioned that
we only use the direct method for solving the resultant linear systems of FIDS (\ref{eq3.15}) and DIDS (\ref{eq3.16}),
respectively, because the maximal size of such resultant linear systems is still \textit{smaller} than 128 and
the superiority of Krylov subspace solvers with circulant preconditioners are slightly less remarkable compared
to the direct solvers in terms of the elapsed CPU time; see e.g., to \cite{Pang12,Gu2017} and the context in the
next example for a discussion.
%
\begin{table}[!htpb]\footnotesize\tabcolsep = 4.8pt
\caption{The $L_{\infty}$- and $L_2$-norm of errors, temporal convergence orders for solving Example 1
with $\alpha = 1.5$, $N(M) = \lceil 2M^{\min\{r\gamma,2-\gamma\}/2} \rceil$, $\mu = 1 +
\frac{\alpha}{2}$, $\kappa_{\mu} = 1$, and $s = 3$.}
\centering
\begin{tabular}{lrrcccrccccr}
\hline &&\multicolumn{5}{c}{DIDS (\ref{eq3.16})} &\multicolumn{5}{c}{FIDS (\ref{eq3.15})}\\
[-2pt]\cmidrule(r{0.7em}r{0.7em}){3-7}\cmidrule(r{0.7em}r{0.6em}){8-12}\\[-11pt]
$(r,\gamma)$ &$M$      &Err$_{\infty}$ &Rate  &Err$_2$  &Rate  &CPU(s)         &Err$_{\infty}$ &Rate  &Err$_2$  &Rate  &CPU(s) \\
\hline
(1,0.8)      &$2^8$    &6.734e-3       &--    &5.974e-3 &--    &0.068          &6.734e-3       &--    &5.974e-3 &--    &\textbf{0.051} \\
             &$2^9$    &3.639e-3       &0.888 &3.173e-3 &0.913 &0.182          &3.639e-3       &0.888 &3.173e-3 &0.913 &\textbf{0.103} \\
             &$2^{10}$ &1.972e-3       &0.884 &1.698e-3 &0.902 &0.619          &1.972e-3       &0.884 &1.698e-3 &0.902 &\textbf{0.204} \\
             &$2^{11}$ &1.108e-3       &0.832 &9.450e-4 &0.845 &2.712          &1.108e-3       &0.832 &9.450e-4 &0.845 &\textbf{0.702} \\
\hline
(2,0.5)      &$2^7$    &4.363e-3       &--    &3.827e-3 &--    &\textbf{0.026} &4.363e-3       &--    &3.827e-3 &--    &0.028          \\
             &$2^8$    &1.964e-3       &1.152 &1.692e-3 &1.177 &\textbf{0.047} &1.964e-3       &1.152 &1.692e-3 &1.177 &0.059          \\
             &$2^9$    &9.497e-4       &1.048 &8.124e-4 &1.059 &0.163          &9.497e-4       &1.048 &8.124e-4 &1.059 &\textbf{0.152} \\
             &$2^{10}$ &4.542e-4       &1.064 &3.934e-4 &1.046 &0.732          &4.542e-4       &1.064 &3.934e-4 &1.046 &\textbf{0.473} \\
\hline
(3,0.8)      &$2^7$    &1.473e-3       &--    &1.304e-3 &--    &\textbf{0.028} &1.473e-3       &--    &1.304e-3 &--    &0.045          \\
             &$2^8$    &5.980e-4       &1.301 &5.251e-4 &1.312 &\textbf{0.074} &5.980e-4       &1.301 &5.251e-4 &1.312 &0.096          \\
             &$2^9$    &2.459e-4       &1.282 &2.136e-4 &1.297 &0.287          &2.459e-4       &1.282 &2.136e-4 &1.297 &\textbf{0.276}          \\
             &$2^{10}$ &1.037e-4       &1.246 &8.954e-5 &1.255 &1.056          &1.037e-4       &1.246 &8.954e-5 &1.255 &\textbf{0.834} \\
\bottomrule
\end{tabular}
\label{tab1}
\end{table}
\begin{table}[t]\footnotesize\tabcolsep = 4.8pt
\caption{The $L_{\infty}$- and $L_2$-norm of errors, temporal convergence orders for solving Example 1
with $\alpha = 0.5$, $N(M) = \lceil 2M^{\min\{r\gamma,2-\gamma\}/2} \rceil$, $\mu = 1 +
\frac{\alpha}{2}$, $\kappa_{\mu} = 1$, and $s = 3$.}
\centering
\begin{tabular}{lrrcccrccccr}
\hline &&\multicolumn{5}{c}{DIDS (\ref{eq3.16})} &\multicolumn{5}{c}{FIDS (\ref{eq3.15})}\\
[-2pt]\cmidrule(r{0.7em}r{0.7em}){3-7}\cmidrule(r{0.7em}r{0.6em}){8-12}\\[-11pt]
$(r,\gamma)$ &$M$      &Err$_{\infty}$ &Rate &Err$_2$   &Rate   &CPU(s)         &Err$_{\infty}$ &Rate  &Err$_2$  &Rate  &CPU(s)          \\
\hline
(1,0.8)      &$2^8$    &2.222e-3       &--    &1.908e-3 &--     &0.058          &2.222e-3       &--     &1.908e-3 &--    &\textbf{0.045} \\
             &$2^9$    &1.240e-3       &0.842 &1.060e-3 &0.847  &0.176          &1.240e-3       &0.842  &1.060e-3 &0.847 &\textbf{0.096} \\
             &$2^{10}$ &6.942e-4       &0.837 &5.918e-4 &0.841  &0.589          &6.942e-4       &0.837  &5.918e-4 &0.841 &\textbf{0.196} \\
             &$2^{11}$ &4.022e-4       &0.787 &3.420e-4 &0.791  &2.654          &4.022e-4       &0.787  &3.420e-4 &0.791 &\textbf{0.558} \\
\hline
(2,0.5)      &$2^7$    &1.608e-3       &--    &1.318e-3 &--     &\textbf{0.023} &1.608e-3       &--     &1.318e-3 &--    &0.032          \\
             &$2^8$    &8.214e-4       &0.969 &6.694e-4 &0.977  &\textbf{0.046} &8.214e-4       &0.969  &6.694e-4 &0.977 &0.061          \\
             &$2^9$    &4.171e-4       &0.978 &3.405e-4 &0.975  &0.181          &4.171e-4       &0.978  &3.405e-4 &0.975 &\textbf{0.147} \\
             &$2^{10}$ &2.112e-4       &0.982 &1.718e-4 &0.987  &0.747          &2.112e-4       &0.982  &1.718e-4 &0.987 &\textbf{0.454} \\
\hline
(3,0.8)      &$2^7$    &4.389e-4       &--    &4.389e-4 &--     &\textbf{0.029} &4.389e-4       &--     &4.389e-4 &--    &0.041          \\
             &$2^8$    &1.869e-4       &1.232 &1.845e-4 &1.250  &\textbf{0.072} &1.869e-4       &1.232  &1.845e-4 &1.250 &0.101          \\
             &$2^9$    &8.001e-5       &1.224 &7.874e-5 &1.228  &0.281          &8.001e-5       &1.224  &7.874e-5 &1.228 &\textbf{0.273} \\
             &$2^{10}$ &3.492e-5       &1.196 &3.439e-5 &1.195  &1.073          &3.492e-5       &1.196  &3.439e-5 &1.195 &\textbf{0.826} \\
\bottomrule
\end{tabular}
\label{tab2}
\end{table}

\begin{table}[!htpb]\footnotesize\tabcolsep = 4.8pt
\caption{The $L_{\infty}$- and $L_2$-norm of errors, spatial convergence orders for solving Example 1
with $\alpha = 1.5$, $M(N) = \lceil (N/2)^{2/\min\{r\gamma,2-\gamma\}} \rceil$, $\mu = 1 +
\frac{\alpha}{2}$, $\kappa_{\mu} = 1$, and $s = 3$.}
\centering
\begin{tabular}{lrrcccrccccr}
\hline &&\multicolumn{5}{c}{DIDS (\ref{eq3.16})} &\multicolumn{5}{c}{FIDS (\ref{eq3.15})}\\
[-2pt]\cmidrule(r{0.7em}r{0.7em}){3-7}\cmidrule(r{0.7em}r{0.6em}){8-12}\\[-11pt]
$(r,\gamma)$ &$N$   &Err$_{\infty}$ &Rate  &Err$_2$  &Rate  &CPU(s)         &Err$_{\infty}$ &Rate  &Err$_2$  &Rate  &CPU(s) \\
\hline
(1,0.8)      &$2^3$ &3.889e-2       &--    &3.802e-2 &--    &0.009          &3.889e-2       &--    &3.802e-2 &--    &\textbf{0.008} \\
             &$2^4$ &8.668e-3       &2.166 &7.760e-3 &2.293 &0.031          &8.668e-3       &2.166 &7.760e-3 &2.293 &\textbf{0.030} \\
             &$2^5$ &1.972e-3	    &2.136 &1.698e-3 &2.193 &0.587          &1.972e-3	    &2.136 &1.698e-3 &2.193 &\textbf{0.256} \\
             &$2^6$ &4.561e-4	    &2.112 &3.846e-4 &2.142	&20.161         &4.561e-4       &2.112 &3.846e-4 &2.142 &\textbf{1.727} \\
\hline
(2,0.5)      &$2^3$ &3.865e-2       &--    &3.789e-2 &--    &\textbf{0.003} &3.865e-2       &--    &3.789e-2 &--    &0.005          \\
             &$2^4$ &8.625e-3       &2.164 &7.734e-3 &2.293 &\textbf{0.010} &8.625e-3       &2.164 &7.734e-3 &2.293 &0.016          \\
             &$2^5$ &1.964e-3       &2.135 &1.692e-3 &2.192 &\textbf{0.045} &1.964e-3       &2.135 &1.692e-3 &2.192 &0.068          \\
             &$2^6$ &4.542e-4       &2.112 &3.834e-4 &2.142 &0.736          &4.542e-4       &2.112 &3.834e-4 &2.142 &\textbf{0.457} \\
\hline
(3,0.8)      &$2^3$ &3.766e-2       &--    &3.791e-2 &--    &\textbf{0.002} &3.766e-2       &--    &3.791e-2 &--    &0.004          \\
             &$2^4$ &8.349e-3       &2.174 &7.710e-3 &2.298	&\textbf{0.006} &8.349e-3       &2.174 &7.710e-3 &2.298	&0.009          \\
             &$2^5$ &1.889e-3	    &2.144 &1.682e-3 &2.197	&\textbf{0.015} &1.889e-3	    &2.144 &1.682e-3 &2.197	&0.032          \\
             &$2^6$ &4.351e-4	    &2.119 &3.801e-4 &2.145	&\textbf{0.108} &4.351e-4	    &2.119 &3.801e-4 &2.145	&0.152          \\
\bottomrule
\end{tabular}
\label{tab3}
\end{table}
\begin{table}[!htpb]\footnotesize\tabcolsep = 4.8pt
\caption{The $L_{\infty}$- and $L_2$-norm of errors, spatial convergence orders for solving Example 1
with $\alpha = 0.5$, $M(N) = \lceil (N/2)^{2/\min\{r\gamma,2-\gamma\}} \rceil$, $\mu = 1 +
\frac{\alpha}{2}$, $\kappa_{\mu} = 1$, and $s = 3$.}
\centering
\begin{tabular}{lrrcccrccccr}
\hline &&\multicolumn{5}{c}{DIDS (\ref{eq3.16})} &\multicolumn{5}{c}{FIDS (\ref{eq3.15})}\\
[-2pt]\cmidrule(r{0.7em}r{0.7em}){3-7}\cmidrule(r{0.7em}r{0.6em}){8-12}\\[-11pt]
$(r,\gamma)$ &$N$   &Err$_{\infty}$ &Rate  &Err$_2$  &Rate  &CPU(s)         &Err$_{\infty}$ &Rate  &Err$_2$  &Rate  &CPU(s)         \\
\hline
(1,0.8)      &$2^3$ &1.197e-2       &--    &1.039e-2 &-- 	&0.007          &1.197e-2       &--    &1.039e-2 &-- 	&\textbf{0.005} \\
             &$2^4$ &2.822e-3	    &2.085 &2.429e-3 &2.097	&0.030          &2.822e-3       &2.085 &2.429e-3 &2.097 &\textbf{0.029} \\
             &$2^5$ &6.942e-4	    &2.023 &5.918e-4 &2.037 &0.583          &6.942e-4       &2.023 &5.918e-4 &2.037 &\textbf{0.247} \\
             &$2^6$ &1.730e-4	    &2.005 &1.467e-4 &2.013	&20.156         &1.730e-4	    &2.005 &1.467e-4 &2.013 &\textbf{1.731} \\
\hline
(2,0.5)      &$2^3$ &1.167e-2	    &--    &1.021e-2 &	    &\textbf{0.003} &1.167e-2       &--    &1.021e-2 &--    &0.004          \\
             &$2^4$ &3.092e-3       &1.917 &2.559e-3 &1.996	&\textbf{0.012} &3.092e-3       &1.917 &2.559e-3 &1.996	&0.014          \\
             &$2^5$ &8.214e-4	    &1.912 &6.694e-4 &1.935	&\textbf{0.044} &8.214e-4	    &1.912 &6.694e-4 &1.935	&0.062          \\
             &$2^6$ &2.112e-4	    &1.960 &1.718e-4 &1.963	&0.734          &2.112e-4	    &1.960 &1.718e-4 &1.963	&\textbf{0.448} \\
\hline
(3,0.8)      &$2^3$ &9.872e-3	    &--    &9.697e-3 &--	&\textbf{0.002} &9.872e-3	    &--    &9.697e-3 &--	&0.004          \\
             &$2^4$ &2.303e-3	    &2.100 &2.264e-3 &2.099	&\textbf{0.005} &2.303e-3	    &2.100 &2.264e-3 &2.099	&0.007          \\
             &$2^5$ &5.584e-4	    &2.044 &5.493e-4 &2.043	&\textbf{0.016} &5.584e-4	    &2.044 &5.493e-4 &2.043	&0.029          \\
             &$2^6$ &1.379e-4	    &2.017 &1.359e-4 &2.016	&\textbf{0.104} &1.379e-4	    &2.017 &1.359e-4 &2.016	&0.139          \\
\bottomrule
\end{tabular}
\label{tab4}
\end{table}

Tables \ref{tab1}--\ref{tab4} present the numerical errors, CPU time (in seconds) and
spatial/temproal convergence rates of both FIDS and DIDS for solving the problem (\ref{eq1.1}),
which satisfy the smooth condition mentioned in Remark \ref{rem3.1}. When we refine the
discretized grid size, it is easily seen that for the temporal direction, the numerical
convergence order is consistent with the theoretical estimate $\mathcal{O}(M^{-\min\{r\gamma,
2-\gamma\}})$ for different $\alpha$'s. Meanwhile, it can find that the numerically spatial
convergence order is exactly consistent with the theoretical estimation $\mathcal{O}(h^2)
$ for different orders of the IFL. In addition, the results of CPU time demonstrate
that the FIDS outperforms the DIDS, especially when the integer $M$ is
increasingly large.
\begin{table}[t]\footnotesize\tabcolsep = 4.8pt
\caption{The $L_{\infty}$- and $L_2$-norm of errors, temporal convergence orders for solving Example 1
with $\alpha = 1.6$, $N(M) = \lceil 2M^{\min\{r\gamma,2-\gamma\}/2} \rceil$, $\mu = 2$, $\kappa_{\mu} = 2$, and $s = 3$}
\centering
\begin{tabular}{lrrcccrccccr}
\hline &&\multicolumn{5}{c}{DIDS (\ref{eq3.16})} &\multicolumn{5}{c}{FIDS (\ref{eq3.15})}\\
[-2pt]\cmidrule(r{0.7em}r{0.7em}){3-7}\cmidrule(r{0.7em}r{0.6em}){8-12}\\[-11pt]
$(r,\gamma)$ &$M$      &Err$_{\infty}$ &Rate  &Err$_2$  &Rate  &CPU(s)         &Err$_{\infty}$ &Rate  &Err$_2$  &Rate  &CPU(s)         \\
\hline
(1,0.8)      &$2^8$    &6.826e-3	   &--    &6.184e-3	&--    &0.067          &6.826e-3	   &--    &6.184e-3	&--    &\textbf{0.048} \\
             &$2^9$    &3.640e-3	   &0.907 &3.238e-3	&0.934 &0.183          &3.640e-3	   &0.907 &3.238e-3	&0.934 &\textbf{0.103} \\
             &$2^{10}$ &1.943e-3	   &0.906 &1.705e-3	&0.925 &0.614          &1.943e-3	   &0.906 &1.705e-3	&0.925 &\textbf{0.199} \\
             &$2^{11}$ &1.075e-3	   &0.854 &9.345e-4	&0.868 &2.710          &1.075e-3	   &0.854 &9.345e-4	&0.868 &\textbf{0.687} \\
\hline
(2,0.5)      &$2^7$    &4.385e-3	   &--    &3.924e-3	&--    &\textbf{0.022} &4.385e-3	   &--    &3.924e-3	&--    &0.028          \\
             &$2^8$    &1.935e-3	   &1.180 &1.701e-3	&1.206 &\textbf{0.044} &1.935e-3	   &1.180 &1.701e-3	&1.206 &0.064          \\
             &$2^9$    &9.178e-4	   &1.076 &8.005e-4	&1.087 &0.176          &9.178e-4	   &1.076 &8.005e-4	&1.087 &\textbf{0.155} \\
             &$2^{10}$ &4.298e-4	   &1.095 &3.697e-4	&1.115 &0.734          &4.298e-4	   &1.095 &3.697e-4	&1.115 &\textbf{0.438} \\
\hline
(3,0.8)      &$2^7$    &1.445e-3	   &--	  &1.304e-3	&--	   &\textbf{0.028} &1.445e-3	   &--	  &1.304e-3	&--	   &0.047          \\
             &$2^8$    &5.716e-4	   &1.338 &5.120e-4	&1.349 &\textbf{0.076} &5.716e-4	   &1.338 &5.120e-4	&1.349 &0.099          \\
             &$2^9$    &2.288e-4	   &1.321 &2.029e-4	&1.336 &0.279          &2.288e-4	   &1.321 &2.029e-4	&1.336 &\textbf{0.271} \\
             &$2^{10}$ &9.389e-5	   &1.285 &8.279e-5	&1.293 &1.059          &9.389e-5	   &1.285 &8.279e-5	&1.293 &\textbf{0.832} \\
\bottomrule
\end{tabular}
\label{tab5}
\end{table}

\begin{table}[!htpb]\footnotesize\tabcolsep = 4.8pt
\caption{The $L_{\infty}$- and $L_2$-norm of errors, temporal convergence orders for solving Example 1
with $\alpha = 0.4$, $N(M) = \lceil 2M^{\min\{r\gamma,2-\gamma\}/2} \rceil$, $\mu = 2$, $\kappa_{\mu} = 2$, and $s = 3$}
\centering
\begin{tabular}{lrrcccrccccr}
\hline &&\multicolumn{5}{c}{DIDS (\ref{eq3.16})} &\multicolumn{5}{c}{FIDS (\ref{eq3.15})}\\
[-2pt]\cmidrule(r{0.7em}r{0.7em}){3-7}\cmidrule(r{0.7em}r{0.6em}){8-12}\\[-11pt]
$(r,\gamma)$ &$M$      &Err$_{\infty}$ &Rate  &Err$_2$  &Rate  &CPU(s)         &Err$_{\infty}$ &Rate  &Err$_2$  &Rate  &CPU(s)        \\
\hline
(1,0.8)      &$2^8$    &1.829e-3	   &--    &1.468e-3	&--    &0.061          &1.829e-3	   &--    &1.468e-3	&--    &\textbf{0.046} \\
             &$2^9$    &1.070e-3	   &0.774 &8.633e-4	&0.766 &0.179          &1.070e-3	   &0.774 &8.633e-4	&0.766 &\textbf{0.094} \\
             &$2^{10}$ &6.206e-4	   &0.785 &5.027e-4	&0.780 &0.592          &6.206e-4	   &0.785 &5.027e-4	&0.780 &\textbf{0.199} \\
             &$2^{11}$ &3.589e-4	   &0.790 &2.911e-4	&0.788 &2.667          &3.589e-4	   &0.790 &2.911e-4	&0.788 &\textbf{0.553} \\
\hline
(2,0.5)      &$2^7$    &1.657e-3	   &--	  &1.343e-3	&--	   &\textbf{0.022} &1.657e-3	   &--	  &1.343e-3	&--	   &0.033          \\
             &$2^8$    &8.409e-4	   &0.978 &6.829e-4	&0.976 &\textbf{0.045} &8.409e-4	   &0.978 &6.829e-4	&0.976 &0.062          \\
             &$2^9$    &4.233e-4	   &0.990 &3.446e-4	&0.987 &0.183          &4.233e-4	   &0.990 &3.446e-4	&0.987 &\textbf{0.151} \\
             &$2^{10}$ &2.128e-4	   &0.992 &1.731e-4	&0.993 &0.744          &2.128e-4	   &0.992 &1.731e-4	&0.993 &\textbf{0.448} \\
\hline
(3,0.8)      &$2^7$    &4.758e-4	   &--    &3.854e-4	&--    &\textbf{0.029} &4.758e-4	   &--    &3.854e-4	&--    &0.042          \\
             &$2^8$    &2.071e-4	   &1.200 &1.687e-4	&1.193 &0.070          &2.071e-4	   &1.200 &1.687e-4	&1.193 &0.101          \\
             &$2^9$    &8.963e-5	   &1.208 &7.315e-5	&1.205 &0.279          &8.963e-5	   &1.208 &7.315e-5	&1.205 &\textbf{0.268} \\
             &$2^{10}$ &3.930e-5	   &1.190 &3.212e-5	&1.187 &1.069          &3.930e-5	   &1.190 &3.212e-5	&1.187 &\textbf{0.821} \\
\bottomrule
\end{tabular}
\label{tab6}
\end{table}

On the other hand, there is another splitting parameter $\mu = 2$ for discretizing the IFL, and
it makes the spatial discretization of IFL enjoy the second-order accuracy \cite{Duo2018a}. According
to Tables \ref{tab5}--\ref{tab6} and Figs. \ref{fig1}--\ref{fig2}, it is not hard to find that
both FIDS and DIDS under such a spatial discretization for solving Example 1 can reach still the
spatial convergence order $\mathcal{O}(h^2)$ and $\mathcal{O}(M^{-\min\{r\gamma,2-\gamma\}})$ with
different settings. Thus such results are also consistent with the theoretical estimate described
in Section \ref{sec3.2}. Again, the results of the elapsed CPU time show that the FIDS
consumes much less time than the DIDS, especially when the integer $M$ is very large.

\vspace{1mm}
\noindent\textbf{Example 2}. The second example is similar to the setting in Example 1, while we choose
$s = 1$ and $\kappa(x,t) = 7[\ln(5 +2x + t) + \cos(xt)]/4$. Then the exact solution $u(x,t)$ and source
term $f(x,t)$ can be computed via the form described in Example 1, and the corresponding initial-boundary
conditions are also similarly obtained.
\begin{table}[!htpb]\footnotesize\tabcolsep = 4.8pt
\caption{The $L_{\infty}$- and $L_2$-norm of errors, temporal convergence orders for solving Example 2
with $\alpha = 1.9$, $N(M) = \lceil 2M^{\min\{r\gamma,2-\gamma\}/\mu} \rceil$, $\mu = 1 + \frac{\alpha}
{2}$ and $\kappa_{\mu} = 1$.}
\centering
\begin{tabular}{lrrcccrccccr}
\hline &&\multicolumn{5}{c}{DIDS (\ref{eq3.16})} &\multicolumn{5}{c}{FIDS (\ref{eq3.15})}\\
[-2pt]\cmidrule(r{0.7em}r{0.7em}){3-7}\cmidrule(r{0.7em}r{0.6em}){8-12}\\[-11pt]
$(r,\gamma)$ &$M$      &Err$_{\infty}$ &Rate  &Err$_2$  &Rate  &CPU(s)         &Err$_{\infty}$ &Rate  &Err$_2$  &Rate  &CPU(s) \\
\hline
(1,0.8)      &$2^8$    &1.724e-2	   &--	  &1.839e-2	&--	   &0.054          &1.724e-2	   &--	  &1.839e-2	&--	   &\textbf{0.046} \\
             &$2^9$    &9.887e-3	   &0.803 &1.056e-2	&0.801 &0.178          &9.887e-3	   &0.803 &1.056e-2	&0.801 &\textbf{0.102} \\
             &$2^{10}$ &5.300e-3	   &0.900 &5.664e-3	&0.898 &0.703          &5.300e-3	   &0.900 &5.664e-3	&0.898 &\textbf{0.256} \\
             &$2^{11}$ &2.997e-3	   &0.823 &3.210e-3	&0.819 &2.586          &2.997e-3	   &0.823 &3.210e-3	&0.819 &\textbf{0.487} \\
\hline
(2,0.5)      &$2^7$    &1.061e-2	   &--	  &1.132e-2	&--	   &\textbf{0.014} &1.061e-2	   &--	  &1.132e-2	&--	   &0.026          \\
             &$2^8$    &5.231e-3	   &1.021 &5.595e-3	&1.017 &\textbf{0.054} &5.231e-3	   &1.021 &5.595e-3	&1.017 &0.059          \\
             &$2^9$    &2.488e-3	   &1.072 &2.668e-3	&1.068 &0.152          &2.488e-3	   &1.072 &2.668e-3	&1.068 &\textbf{0.141} \\
             &$2^{10}$ &1.241e-3	   &1.004 &1.334e-3	&1.001 &0.721          &1.241e-3	   &1.004 &1.334e-3	&1.001 &\textbf{0.457} \\
\hline
(3,0.8)      &$2^7$    &3.951e-3	   &--	  &4.234e-3	&--	   &\textbf{0.019} &3.951e-3	   &--	  &4.234e-3	&--	   &0.046          \\
             &$2^8$    &1.645e-3	   &1.264 &1.767e-3	&1.260 &\textbf{0.066} &1.645e-3	   &1.264 &1.767e-3	&1.260 &0.098          \\
             &$2^9$    &6.895e-4	   &1.255 &7.427e-4	&1.251 &0.272          &6.895e-4	   &1.255 &7.427e-4	&1.251 &\textbf{0.267} \\
             &$2^{10}$ &2.851e-4	   &1.274 &3.079e-4	&1.270 &1.100          &2.851e-4	   &1.274 &3.079e-4	&1.270 &\textbf{0.824} \\
\bottomrule
\end{tabular}
\label{tab7}
\end{table}
\begin{table}[t]\footnotesize\tabcolsep = 4.8pt
\caption{The $L_{\infty}$- and $L_2$-norm of errors, temporal convergence orders for solving Example 2 with $\alpha = 0.5$,
$N(M) = \lceil 2M^{\min\{r\gamma,2-\gamma\}/\mu} \rceil$, $\mu = 1 + \frac{\alpha}{2}$ and $\kappa_{\mu} = 1$.} 
\centering 
\begin{tabular}{lrrcrccrccccrrrrr}
\hline &&\multicolumn{5}{c}{DIDS (\ref{eq3.16})} &\multicolumn{5}{c}{FIDS (\ref{eq3.15})}\\
[-2pt]\cmidrule(r{0.7em}r{0.7em}){3-7}\cmidrule(r{0.7em}r{0.6em}){8-12}\\[-11pt]
$(r,\gamma)$ &$M$      &Err$_{\infty}$ &Rate  &Err$_2$  &Rate  &CPU(s)         &Err$_{\infty}$ &Rate  &Err$_2$  &Rate  &CPU(s)         \\
\hline
(1,0.8)      &$2^8$    &1.699e-3	   &--	  &1.698e-3	&--	   &0.077          &1.699e-3	   &--	  &1.698e-3	&--	   &0.084          \\
             &$2^9$    &1.016e-3	   &0.742 &1.011e-3	&0.747 &0.321          &1.016e-3	   &0.742 &1.011e-3	&0.747 &\textbf{0.217} \\
             &$2^{10}$ &6.013e-4	   &0.757 &5.980e-4	&0.758 &1.235          &6.013e-4	   &0.757 &5.980e-4	&0.758 &\textbf{0.660} \\
             &$2^{11}$ &3.521e-4	   &0.772 &3.496e-4	&0.774 &5.899          &3.521e-4	   &0.772 &3.496e-4	&0.774 &\textbf{2.759} \\
\hline
(2,0.5)      &$2^7$    &1.619e-3	   &--	  &1.609e-3	&--	   &\textbf{0.031} &1.619e-3	   &--	  &1.609e-3	&--	   &0.050          \\
             &$2^8$    &8.238e-4	   &0.975 &8.174e-4	&0.977 &\textbf{0.158} &8.238e-4	   &0.975 &8.174e-4	&0.977 &0.173          \\
             &$2^9$    &4.189e-4	   &0.976 &4.157e-4	&0.975 &0.979          &4.189e-4	   &0.976 &4.157e-4	&0.975 &\textbf{0.887} \\
             &$2^{10}$ &2.115e-4	   &0.986 &2.098e-4	&0.987 &10.409         &2.115e-4	   &0.986 &2.098e-4	&0.987 &\textbf{9.215} \\
\hline
(3,0.8)      &$2^5$    &1.274e-3	   &--	  &1.340e-3	&--	   &\textbf{0.011} &1.274e-3	   &--	  &1.340e-3	&--	   &0.012          \\
             &$2^6$    &5.719e-4	   &1.156 &5.970e-4	&1.166 &\textbf{0.025} &5.719e-4	   &1.156 &5.970e-4	&1.166 &0.038          \\
             &$2^7$    &2.529e-4	   &1.178 &2.627e-4	&1.185 &\textbf{0.102} &2.529e-4	   &1.178 &2.627e-4	&1.185 &0.130          \\
             &$2^8$    &1.109e-4	   &1.189 &1.149e-4	&1.193 &1.426          &1.109e-4	   &1.189 &1.149e-4	&1.193 &\textbf{1.347} \\
\bottomrule
\end{tabular}
\label{tab8}
\end{table}
\begin{table}[t]\footnotesize\tabcolsep = 4.8pt
\caption{The $L_{\infty}$- and $L_2$-norm of errors, spatial convergence orders  for solving Example 2 with $\alpha = 1.9$,
$M(N) = \lceil (N/2)^{\mu/\min\{r\gamma,2-\gamma\}} \rceil$, $\mu = 1 + \frac{\alpha}{2}$ and $\kappa_{\mu} = 1$.}
\centering
\begin{tabular}{lrrcccrccccrrrrr} 
\hline &&\multicolumn{5}{c}{DIDS (\ref{eq3.16})} &\multicolumn{5}{c}{FIDS (\ref{eq3.15})}\\
[-2pt]\cmidrule(r{0.7em}r{0.7em}){3-7}\cmidrule(r{0.7em}r{0.6em}){8-12}\\[-11pt]
$(r,\gamma)$ &$N$   &Err$_{\infty}$ &Rate  &Err$_2$  &Rate  &CPU(s)         &Err$_{\infty}$ &Rate  &Err$_2$  &Rate  &CPU(s)         \\
\hline
(1,0.8)      &9     &7.770e-2	    &--	   &8.304e-2 &--	&\textbf{0.006} &7.770e-2	    &--	   &8.304e-2 &--	&0.007          \\
             &18    &1.929e-2	    &2.010 &2.052e-2 &2.017	&0.035          &1.929e-2	    &2.010 &2.052e-2 &2.017	&\textbf{0.033} \\
             &36    &4.719e-3	    &2.031 &5.045e-3 &2.024	&0.834          &4.719e-3	    &2.031 &5.045e-3 &2.024	&\textbf{0.254} \\
             &72    &1.153e-3	    &2.033 &1.238e-3 &2.026	&23.551         &1.153e-3	    &2.033 &1.238e-3 &2.026	&\textbf{1.927} \\
\hline
(2,0.5)      &9     &7.663e-2	    &--	   &8.193e-2 &-- 	&\textbf{0.002} &7.663e-2	    &--	   &8.193e-2 &-- 	&0.005          \\
             &18    &1.903e-2	    &2.010 &2.026e-2 &2.016	&\textbf{0.008} &1.903e-2	    &2.010 &2.026e-2 &2.016	&0.014          \\
             &36    &4.657e-3	    &2.031 &4.983e-3 &2.023	&\textbf{0.053} &4.657e-3	    &2.031 &4.983e-3 &2.023	&0.077          \\
             &72    &1.138e-3	    &2.033 &1.223e-3 &2.026 &0.879          &1.138e-3	    &2.033 &1.223e-3 &2.026 &\textbf{0.426} \\
\hline
(3,0.8)      &9     &7.663e-2	    &--	   &8.195e-2 &--	&\textbf{0.002} &7.663e-2	    &--	   &8.195e-2 &--	&0.004          \\
             &18    &1.902e-2	    &2.011 &2.025e-2 &2.017	&\textbf{0.005} &1.902e-2	    &2.011 &2.025e-2 &2.017	&0.008          \\
             &36    &4.650e-3	    &2.032 &4.977e-3 &2.025	&\textbf{0.016} &4.650e-3	    &2.032 &4.977e-3 &2.025	&0.033          \\
             &72    &1.135e-3	    &2.034 &1.221e-3 &2.028	&\textbf{0.114} &1.135e-3	    &2.034 &1.221e-3 &2.028	&0.162          \\
\bottomrule
\end{tabular}
\label{tab9}
\end{table}
\begin{table}[!htpb]\footnotesize\tabcolsep = 4.8pt
\caption{The $L_{\infty}$- and $L_2$-norm of errors, spatial convergence orders  for solving Example 2 with $\alpha = 0.5$,
$M(N) = \lceil (N/2)^{\mu/\min\{r\gamma,2-\gamma\}} \rceil$, $\mu = 1 + \frac{\alpha}{2}$ and $\kappa_{\mu} = 1$}
\centering
\begin{tabular}{lrrcccrccccrrrrr} 
\hline &&\multicolumn{5}{c}{DIDS (\ref{eq3.16})} &\multicolumn{5}{c}{FIDS (\ref{eq3.15})}\\
[-2pt]\cmidrule(r{0.7em}r{0.7em}){3-7}\cmidrule(r{0.7em}r{0.6em}){8-12}\\[-11pt]
$(r,\gamma)$ &$N$      &Err$_{\infty}$ &Rate  &Err$_2$  &Rate  &CPU(s) &Err$_{\infty}$ &Rate  &Err$_2$  &Rate  &CPU(s)         \\
\hline
(1,0.8)      &$2^6$    &1.871e-3	   &--	  &1.871e-3	&--    &0.065          &1.871e-3	   &--	  &1.871e-3	&--    &\textbf{0.063} \\
             &$2^7$    &8.355e-4	   &1.163 &8.323e-4	&1.169 &0.539          &8.355e-4	   &1.163 &8.323e-4	&1.169 &\textbf{0.364} \\
             &$2^8$    &3.642e-4	   &1.198 &3.616e-4	&1.203 &5.365          &3.642e-4	   &1.198 &3.616e-4	&1.203 &\textbf{2.438} \\
             &$2^9$    &1.556e-4	   &1.227 &1.543e-4	&1.229 &88.190         &1.556e-4	   &1.227 &1.543e-4	&1.229 &\textbf{51.194} \\
\hline
(2,0.5)      &$2^6$    &2.666e-3	   &--	  &2.654e-3	&--	   &\textbf{0.013} &2.666e-3	   &--	  &2.654e-3	&--	   &0.023          \\
             &$2^7$    &1.157e-3	   &1.204 &1.149e-3	&1.208 &\textbf{0.073} &1.157e-3	   &1.204 &1.149e-3	&1.208 &0.091          \\
             &$2^8$    &4.970e-4	   &1.219 &4.933e-4	&1.219 &0.584          &4.970e-4	   &1.219 &4.933e-4	&1.219 &\textbf{0.571} \\
             &$2^9$    &2.115e-4	   &1.232 &2.098e-4	&1.234 &9.855          &2.115e-4	   &1.232 &2.098e-4	&1.234 &\textbf{9.277} \\
\hline
(3,0.8)      &$2^5$    &2.581e-3	   &--	  &2.739e-3	&--	   &\textbf{0.003} &2.581e-3	   &--	  &2.739e-3	&--	   &0.006          \\
             &$2^6$    &1.114e-3	   &1.211 &1.169e-3	&1.228 &\textbf{0.007} &1.114e-3	   &1.211 &1.169e-3	&1.228 &0.014          \\
             &$2^7$    &4.677e-4	   &1.253 &4.874e-4	&1.262 &\textbf{0.025} &4.677e-4	   &1.253 &4.874e-4	&1.262 &0.046          \\
             &$2^8$    &2.000e-4	   &1.226 &2.075e-4	&1.232 &\textbf{0.161} &2.000e-4	   &1.226 &2.075e-4	&1.232 &0.232          \\
\bottomrule
\end{tabular}
\label{tab10}
\end{table}

In this example, it notes that the exact solution $u(x,\cdot)\in\mathcal{C}^{1,\frac{\alpha}{2}}(\mathbb{R})$ satisfies
the less smoother condition than that in Example 1. Moreover, it is seen
from Tables \ref{tab7}--\ref{tab10} that for the temporal direction, the
numerical convergence rate of both FIDS and DIDS is consistent with the theoretical estimate
$\mathcal{O}(M^{-\min\{r\gamma,2 - \gamma\}})$ for different settings. However, it remarks that the spatial convergence
rate of both FIDS and DIDS can at least approach to $1 + \frac{\alpha}{2}$, especially when $\alpha$ increasingly goes
to 2, and the spatial convergence orders of both FIDS and DIDS are almost 2. These results on spatial convergence rate of
both FIDS and DIDS are fairly better than the theoretical estimate in Remark \ref{rem3.1}. It implies that the error
analysis and smooth condition of the numerical discretization of IFL used to establish the IDS can be further sharped
and weakened, respectively. Analogously, the average CPU time of FIDS is smaller than that of DIDS for
problem (\ref{eq1.1}), when the number of time levels is increasingly large.

\begin{figure}[!htpb]
\centering
\subfigure[The number of time levels $M$ versus the number of spatial grid nodes ($N$) for
$(r,\gamma) = (1,0.8)$.]{
\label{fig:subfig:a} 
\includegraphics[width=3.15in]{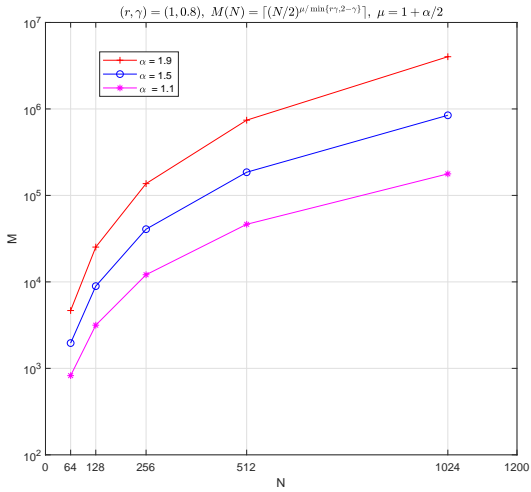}}
\subfigure[The eigenvalue distributions of original and preconditioned matrices when $\alpha = 1.9$.]{
\label{subfig:b} 
\includegraphics[width=3.15in]{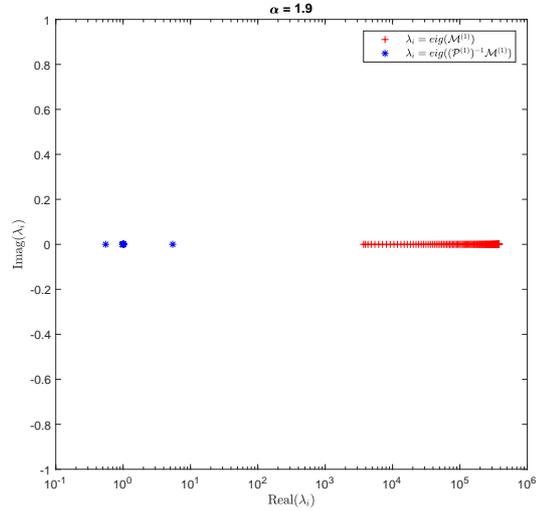}}\\
\subfigure[The eigenvalue distributions of original and preconditioned matrices when $\alpha = 1.5$.]{
\label{subfig:c} 
\includegraphics[width=3.15in]{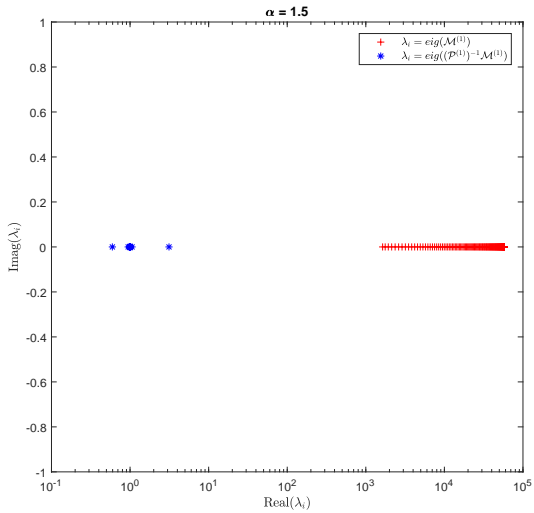}}
\subfigure[The eigenvalue distributions of original and preconditioned matrices when $\alpha = 1.1$.]{
\label{subfig:d} 
\includegraphics[width=3.15in]{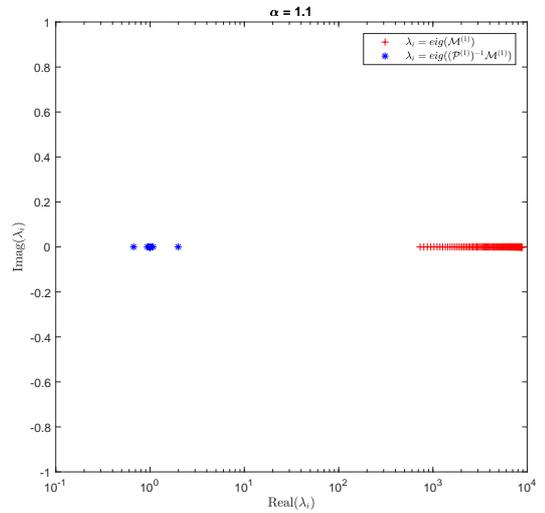}}
\caption{The brief performance analysis of FIDS (\ref{eq3.15}) and DIDS (\ref{eq3.16}) with no/circulant
preconditioners for solving Example 2 with $N = 2^7$ in the subfigures: (b)--(d).}
\label{fig3} 
\end{figure}

\begin{sidewaystable}[!htpb]\footnotesize\tabcolsep = 4.8pt
\caption{The errors and performance of FIDS and DIDS with direct, iterative and preconditioned iterative solvers for
Example 2 with different settings and $M(N) = \lceil (N/2)^{\mu/\min\{r\gamma,2-\gamma\}} \rceil$, $\mu = 1 + \frac{\alpha}{2}$ and $\kappa_{\mu} = 1$.}
\centering
\begin{tabular}{lrrcrrrrrrrrrr}
\hline &&&&\multicolumn{1}{c}{DIDS}&\multicolumn{1}{c}{FIDS} &\multicolumn{2}{c}{DIDS + $\mathcal{I}$}&\multicolumn{2}{c}{DIDS + $\mathcal{P}$}
&\multicolumn{2}{c}{FIDS + $\mathcal{I}$}&\multicolumn{2}{c}{FIDS + $\mathcal{P}$}\\
[-2pt]\cmidrule(r{0.7em}r{0.7em}){5-5}\cmidrule(r{0.7em}r{0.7em}){6-6}\cmidrule(r{0.7em}r{0.6em}){7-8}
\cmidrule(r{0.7em}r{0.6em}){9-10}\cmidrule(r{0.7em}r{0.6em}){11-12}\cmidrule(r{0.7em}r{0.6em}){13-14}\\[-11pt]
$(r,\gamma,\alpha)$ &$N$    &Err$_{\infty}$ &Err$_2$  &CPU(s)          &CPU(s)         &Its   &CPU(s)  &Its &CPU(s)  &Its   &CPU(s) &Its  &CPU(s)                   \\
\hline
(2,0.5,1.5)         &$2^6$    &6.362e-4	    &7.727e-4 &0.162           &\textbf{0.156} &44.1  &0.658   &7.6 &0.332   &44.1  &0.655  &7.7  &0.347           \\
                    &$2^7$    &1.542e-4	    &1.889e-4 &1.793           &\textbf{1.014} &66.8  &5.401   &7.9 &2.380   &66.7  &4.693  &7.9  &1.689           \\
                    &$2^8$    &4.469e-5	    &4.643e-5 &22.217          &7.631          &92.7  &45.111  &8.0 &20.777  &92.7  &30.132 &8.0  &\textbf{6.698}  \\
                    &$2^9$    &1.333e-5	    &1.236e-5 &413.96          &170.24         &125.7 &570.75  &8.2 &285.50  &125.7 &306.26 &8.2  &\textbf{43.338} \\
(3,0.8,1.5)         &$2^7$    &1.381e-4	    &1.734e-4 &\textbf{0.286}  &0.372          &41.5  &0.935   &6.6 &0.446   &41.6  &0.900  &6.6  &0.439            \\
                    &$2^8$    &3.265e-5	    &4.164e-5 &2.931           &2.142          &48.6  &4.479   &6.9 &2.104   &48.6  &4.091  &6.9  &\textbf{1.545}   \\
                    &$2^9$    &7.897e-6	    &9.997e-6 &44.728          &33.845         &54.0  &35.156  &7.0 &16.488  &53.9  &26.938 &7.1  &\textbf{8.143}   \\
                    &$2^{10}$ &2.319e-6	    &2.396e-6 &624.02          &511.16         &57.8  &206.13  &7.0 &138.18  &57.8  &97.137 &7.0  &\textbf{34.543}  \\
\hline
(2,0.5,1.9)         &$2^6$    &1.446e-3	    &1.553e-3 &0.507           &\textbf{0.386} &65.6  &2.123   &7.9 &0.853   &65.6  &1.879  &7.9  &0.748            \\
                    &$2^7$    &3.533e-4	    &3.811e-4 &7.831           &\textbf{2.087} &116.4 &22.072  &7.9 &8.574   &116.4 &17.458 &7.9  &3.838            \\
                    &$2^8$    &8.636e-5	    &9.354e-5 &194.74          &21.665         &198.3 &275.11  &8.1 &135.58  &198.3 &203.71 &8.1  &\textbf{18.073}  \\
                    &$2^9$    &2.112e-5	    &2.297e-5 &48229.2         &1029.78        &316.1 &11093.2 &8.1 &2474.12 &316.3 &2620.56 &8.1  &\textbf{132.04} \\
(3,0.8,1.9)         &$2^7$    &3.521e-4	    &3.800e-4 &0.935           &\textbf{0.689} &73.8  &3.317   &7.3 &1.201   &73.9  &3.222  &7.3  &1.049            \\
                    &$2^8$    &8.596e-5	    &9.318e-5 &10.468          &5.134          &100.3 &20.579  &7.3 &7.464   &100.2 &17.025 &7.3  &\textbf{3.679}   \\
                    &$2^9$    &2.100e-5	    &2.286e-5 &173.492         &84.047         &127.3 &246.74  &7.4 &81.104  &127.5 &151.45 &7.4  &\textbf{21.436}  \\
                    &$2^{10}$ &5.136e-6	    &5.614e-6 &2438.12         &1505.84        &151.0 &1616.83 &7.3 &1003.78 &151.0 &714.72 &7.3  &\textbf{86.540}  \\
\bottomrule
\end{tabular}
\label{tab11}
\end{sidewaystable}
\begin{figure}[!htpb]
\centering
\includegraphics[width=3.15in]{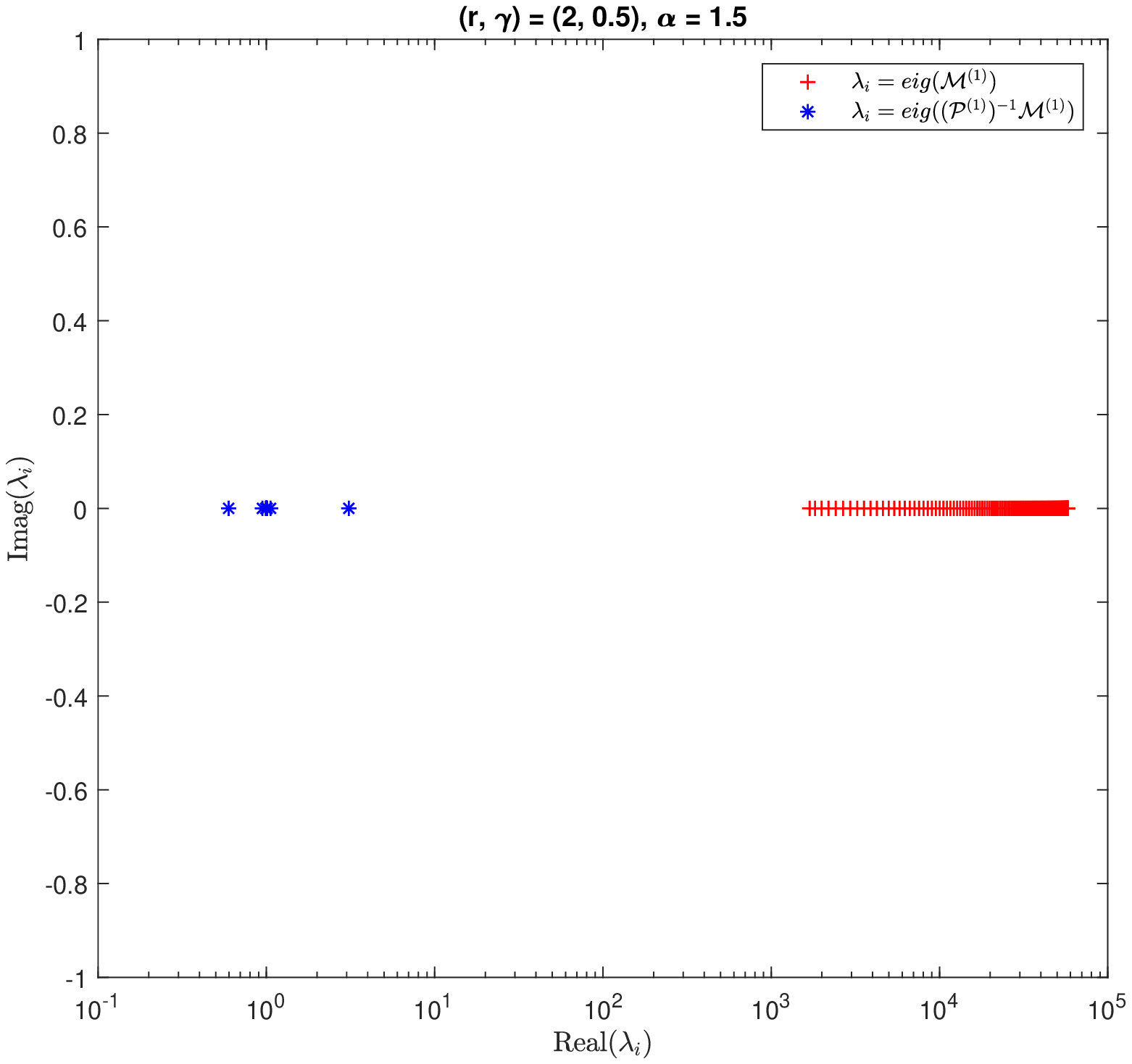}
\includegraphics[width=3.15in]{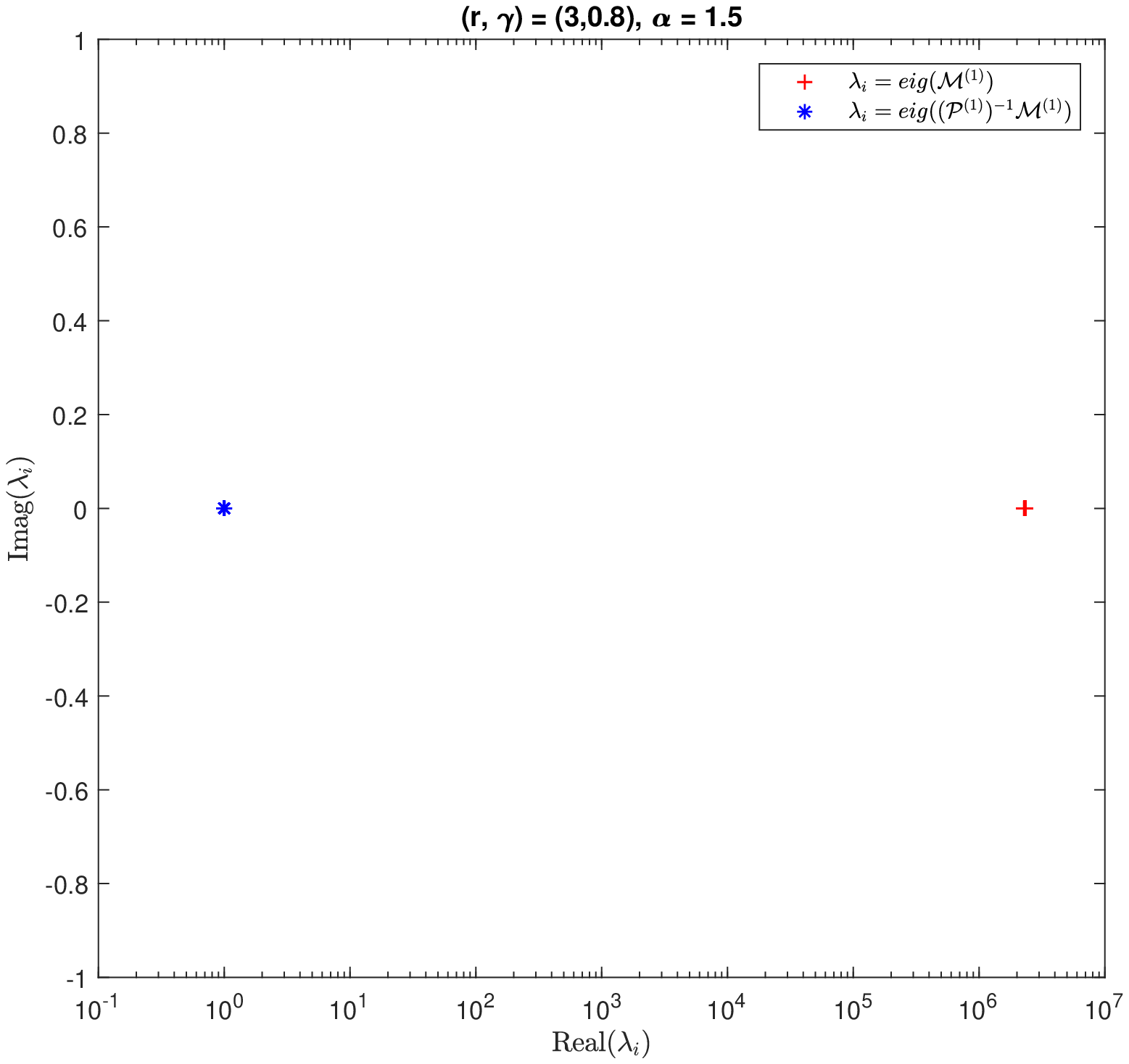}
\includegraphics[width=3.15in]{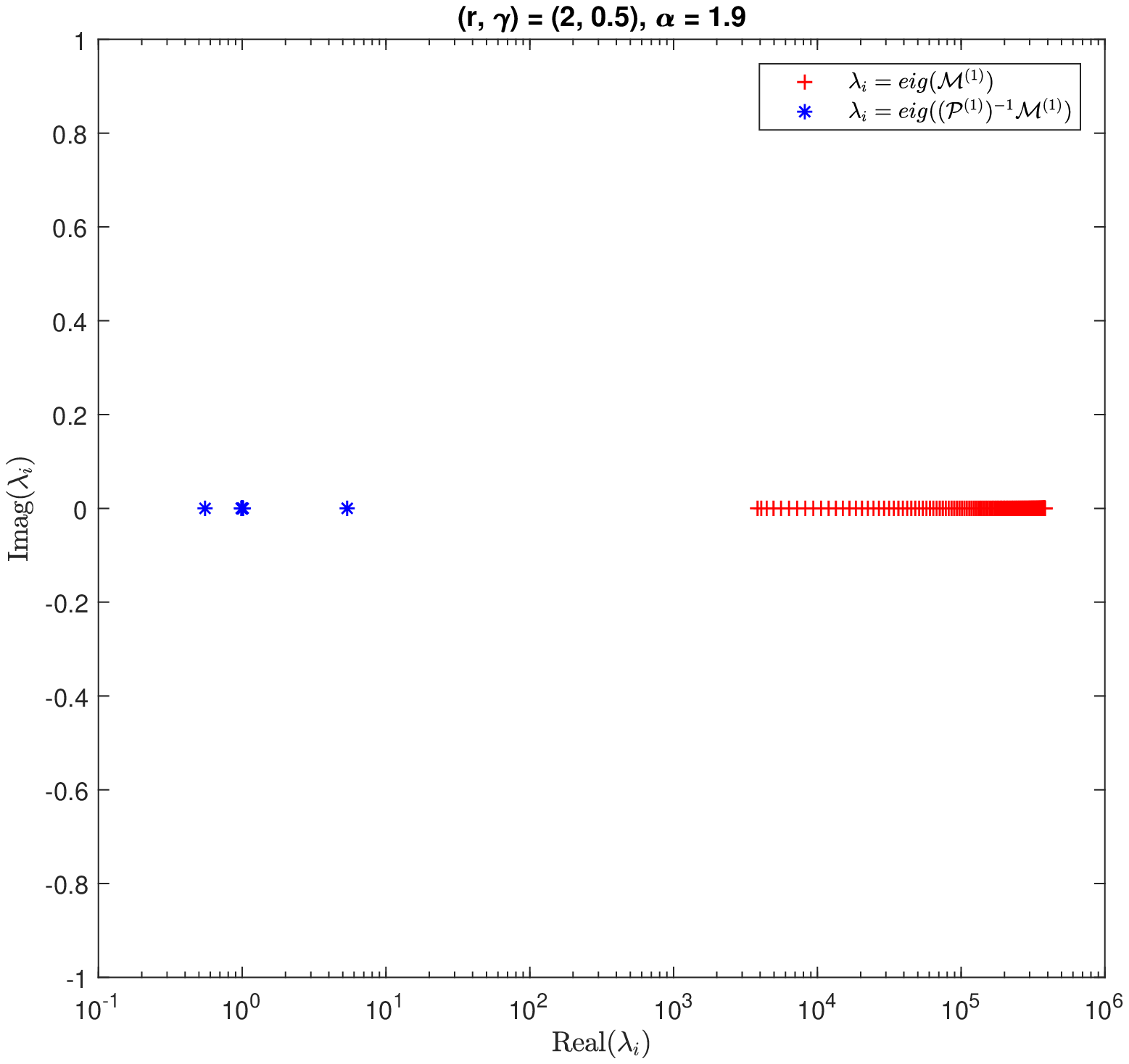}
\includegraphics[width=3.15in]{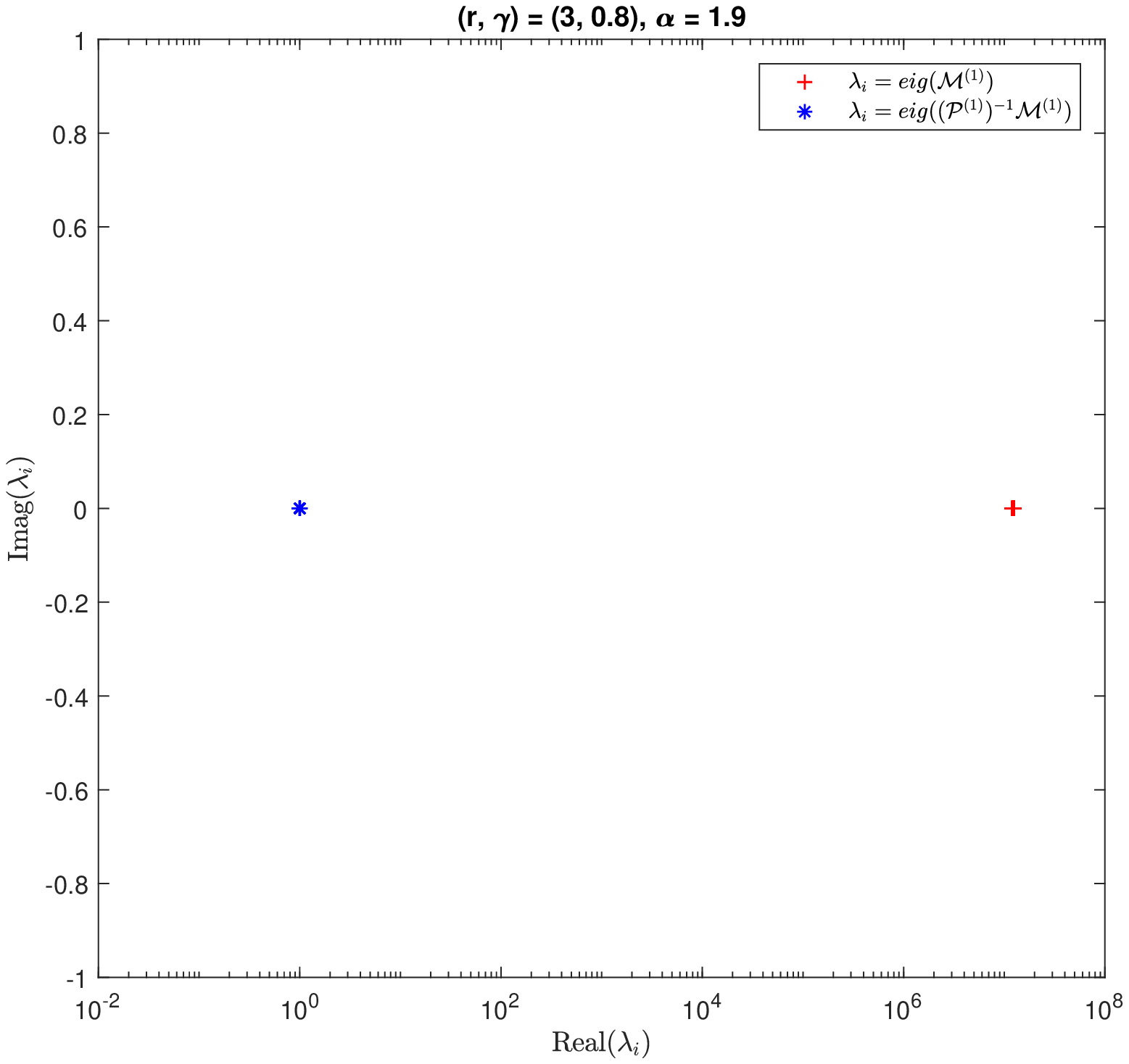}
\caption{Eigenvalue distributions of the original and preconditioned matrices from both FIDS (\ref{eq3.15}) and DIDS (\ref{eq3.16}) for Example 2 with
different setting and $N = 2^7$.}
\label{fig4}
\end{figure}

On the other hand, Table \ref{tab11} and Figs. \ref{fig3}--\ref{fig4} are carried
out to show the effectiveness of the proposed circulant preconditioners,
which is especially useful for the order $\alpha~(\rightarrow 2)$ of IFL; refer to
Figs. \ref{fig3}--\ref{fig4} as well. For Fig. \ref{fig3}(a), it implies that if we increase $N$, then the number of time level will be too huge to make a concise
comparison of FIDS and DIDS (with no/circulant preconditioners). Moreover, due to the large number $M$, the family of
FIDS should be more efficient than the counterparts of DIDS for solving the problem (\ref{eq1.1}). As seen from Table
\ref{tab11}, it finds that the proposed circulant preconditioner is efficient to accelerate the implementations of both
FIDS and DIDS in terms of the reduction of ``Its" and ``CPU", especially for large integers $M$ and $N$. This observation
can be also supported by the clustering eigenvalue distributions shown in Figs. \ref{fig3}--\ref{fig4}. Moreover, the
number of iterations of ``DIDS + $\mathcal{P}$" and ``FIDS + $\mathcal{P}$" is roughly independent of decreasing spatial
grid size. The above results of circulant preconditioners are exactly consistent with the theoretical investigations given
in Section \ref{sec4}. In one word, the ``FIDS + $\mathcal{P}$" is the most promising numerical method for solving the problem
(\ref{eq1.1}), especially with large integers $M,~N$ and $M > N$.
%
\section{Conclusions}
\label{sec6}
In this work, we proposed two fast and easy-to-implement IDSs (i.e., FIDS and DIDS) for solving the
TSFDE (\ref{eq1.1}) with non-smooth initial data, which was not well-stuided in the previous work.
Meanwhile, both the solvability, stability and convergence rate of the proposed IDSs with non-uniform
temproal steps are rigorously proved via the matrix properties, which are meticulously derived from the
direct discretization of IFL. Numerical results in Section \ref{sec5} are reported to support our theoretical
findings. In addition, although the focus is the one-dimensional
spatial domain in this work, we note that the proposed methods utilizing spatial discretizations \cite{Duo2019z}
can be directly adapted and corresponding results remain valid for two- and three-dimensional cases,
which will be precisely presented in our another coming manuscript.

On the other hand, due to the nonlocality of Caputo fractional derivative in the TSFDE (\ref{eq1.1}),
the numerical scheme needs to repeat the weighted sum of solutions {\color{red}at} previous time levels. In order to
reduce the computational cost, we exploit the fast SOE approximation of graded $L1$ formula to result
in the FIDS, which {\color{red}is} cheaper than the DIDS, especially for large integer $M$. However, no
matter FIDS or DIDS, they both need to solve the dense discretized systems, which are still time-consuming.
It implies that the efficient implementation of FIDS and DIDS should be further considered. With the help
of Toeplitz-like matrix, we construct the BiCGSTAB with circulant preconditioners for solving the series
of discretized linear systems (cf. Eq. \ref{eq3.16gu} and Eq. \ref{eq4.1})) without storing any matrices.
It makes the FIDS (or DIDS) only require $\mathcal{O}(NN_{exp})$ (or $\mathcal{O}(NM)$) memory requirement and
$\mathcal{O}(MN(\log N + N_{exp}))$ (or $\mathcal{O}(MN(\log N + M))$) computational complexity. To ensure
circulant preconditioners efficient, we theoretically show that the eigenvalues of preconditioned matrices
cluster around 1, expect for few outliers. The vast majority of these eigenvalues are well separated away
from 0. It means that the BiCGSTAB with circulant preconditioners for solving the discretized linear systems
can converge very fast. Numerical experiments are reported to show the effectiveness of the FIDS and DIDS
with PBiCGSTAB solvers in terms of the elapsed CPU time and number of iterations, especially
the former one.

Finally, it is meaningful to note that the numerically spatial convergence order of FIDS or DIDS is better
than the theoretical estimate of FIDS or DIDS, see Example 2. It means that the error analysis
of numerical schemes for solving time-dependent problems is different from the numerical IFL described in \cite{Duo2018a}, because their result is the error estimate for the IFL, but in current work it is the
estimate for solutions of problem (\ref{eq1.1}). The more refined error/convergence analysis of FIDS or DIDS is worth exploring in
our future work; refer e.g., to
\cite{Duo2019z,Acosta3,Bonito19} for a discussion. Our current work includes applying the FIDS and DIDS for
solving the (nonlinear) multi-dimensional TSFDEs (on unbounded domains) with nonhomogeneous
boundary conditions, and designing the more efficient preconditioning techniques, such as $\tau$-algebra,
multigrid \cite{Pang12} and banded preconditioners, for the corresponding two-
and three-level Toeplitz discrtized linear systems.



\section*{Acknowledgement}
{\em We are very grateful to the anonymous referees for their invaluable comments
and insightful suggestions that have
greatly improved the presentation of this paper. This research is supported by NSFC (11801463), the Applied Basic
Research Project of Sichuan Province (2020YJ0007), the Fundamental Research
Funds for the Central Universities (JBK1902028), the Science and Technology
Development Fund, Macau SAR (file no. 0118/2018/A3), MYRG2018-00015-FST
from University of Macau, and the US National Science Foundation (DMS-1620465).}
\vspace{5mm}

\noindent\textbf{References}

\end{document}